\documentclass[11pt]{article}
\usepackage{geometry}
\usepackage{epsfig}
\usepackage[linesnumbered,ruled,vlined]{algorithm2e}
\usepackage{xcolor}
\usepackage{varwidth}
\usepackage{lineno, hyperref}
\usepackage{graphicx}
\usepackage{epstopdf}
\usepackage{amsmath}
\usepackage{amsthm}
\usepackage{romannum}
\usepackage{amsthm}

\usepackage{booktabs}
\usepackage{float}
\usepackage{subcaption}
\usepackage{caption}
\usepackage{authblk}
\usepackage{abstract}

\date{}
\usepackage{cleveref}
\usepackage[figurename=Fig.]{caption}

\newtheorem{thm}{Theorem}[section]

\newtheorem{example}[thm]{Example}

\numberwithin{equation}{section}

\begin{document}
\providecommand{\keywords}[1]{\textit{Keywords- } #1}

\renewcommand{\thepage}{\arabic{page}}
\setcounter{page}{1}
\title{\textbf{Predator Prey Scavenger Model using Holling's Functional Response of Type III and Physics-Informed Deep Neural Networks}}

\author[1]{Aneesh Panchal\thanks{Email: \texttt{aneeshpanchal840850@gmail.com}}}
\author[2]{Kirti Beniwal\thanks{Email: \texttt{kirtibnwl1912@gmail.com}}}
\author[2]{Vivek Kumar\thanks{Email: \texttt{vivekkumar.ag@gmail.com}}}

\affil[1]{Department of Computational and Data Sciences, Indian Institute of Science, Bangalore, Karnataka$-$560012, INDIA}
\affil[2]{Department of Applied Mathematics, Delhi Technological University, Delhi$-$110042, INDIA}
\maketitle

\begin{abstract}
Nonlinear mathematical models introduce the relation between various physical and biological interactions present in nature. One of the most famous models is the Lotka-Volterra model which defined the interaction between predator and prey species present in nature. However, predators, scavengers, and prey populations coexist in a natural system where scavengers can additionally rely on the dead bodies of predators present in the system. Keeping this in mind, the formulation and simulation of the predator prey scavenger model is introduced in this paper. For the predation response, respective prey species are assumed to have Holling’s functional response of type III. The proposed model is tested for various simulations and is found to be showing satisfactory results in different scenarios. After simulations, the American forest dataset is taken for parameter estimation which imitates the real-world case. For parameter estimation, a physics-informed deep neural network is used with the Adam backpropagation method which prevents the avalanche effect in trainable parameters updation. For neural networks, mean square error and physics-informed informed error are considered. After the neural network, the hence-found parameters are fine-tuned using the Broyden–Fletcher–Goldfarb–Shanno algorithm. Finally, the hence-found parameters using a natural dataset are tested for stability using Jacobian stability analysis. Future research work includes minimization of error induced by parameters, bifurcation analysis, and sensitivity analysis of the parameters.
\end{abstract}

\section{Introduction}\label{introduction}
Mathematical modeling and differential equations increase our understanding of population dynamics and ecology over time. Differential equations \cite{applications} are very useful for defining various natural phenomena such as simple harmonic motion (SHM), Newton’s law of cooling, chemical equations, radioactive decay, and many more. But the study of differential equations increased drastically when the study of environment meets mathematics that is population dynamics \cite{Preyrefugee}, mathematical ecology \cite{Trophicchains}, epidemiology \cite{diseaseprey} etc. Population dynamics \cite{popdynamics1, popdynamics2, popdynamics3} is the study of the variation of the population over time according to different natural parameters. Several continuous variables can be studied jointly to generate a system of ordinary differential equations. It is generally necessary to estimate and fine-tune specific parameters to make the dynamical system valid, providing reliable mathematical explanations for certain phenomena or predictions. In most cases, estimating the parameters involved in dynamical systems is, in fact, a challenging optimization problem that needs special consideration due to some unknown natural parameters. The study of population dynamics got a boost when Italian mathematician \textit{Vito Volterra} \cite{LotkaVolterra} introduced a differential equation model for the population dynamics of two different interacting species viz. predator and prey in nature. Predator-prey modeling using the Lotka-Volterra classical model and Holling's functional response has been a great instrument in advancing the understanding of ecological systems. While the Lotka-Volterra model provides a fundamental basis for studying population interactions, Holling's functional response enhances its applicability by considering the complexities of predator-feeding behavior. Integrating these models contributes to a more comprehensive and realistic representation of predator-prey dynamics in ecological modeling. Further advancements in this field are crucial for addressing contemporary ecological challenges and improving our ability to manage and conserve diverse ecosystems.
\par
In exploring the dynamics of predator-prey models, several key research papers provide valuable insights. One such study by Sarwardi et.al. \cite{Preyrefugee}, dives into the behavior of a two-predator model with prey refugees. By considering the adaptive strategy of prey-seeking refugees, the research provides insights on how refugee behavior affects the stability of the system. Another study by Kant and Kumar \cite{diseaseprey}, focuses on the stability analysis of a predator-prey system where prey populations migrate, and both species are susceptible to disease. This research paper deepens our understanding of migration and disease integration with predator-prey population variations and its impact on the stability of the system. Also, a study by Had{\v{z}}iabdi{\'c} et.al. \cite{LotVol} extends the classical Lotka-Volterra model by adding two predators and common prey, providing insights into the dynamics and stability analysis of multiple predators with a same prey population.
\par
Further research by Wu et.al. \cite{figures} on the delayed predator-prey system, consider factors such as fear effect, herd behavior, and disease in the susceptible prey. The inclusion of time delays and additional ecological aspects offers a more realistic representation of predator-prey interactions, which increases our understanding of the natural behavior of predator and prey. Research by Arora and Kumar \cite{delayholling}, introduces a delayed version of a prey-predator system with a modified Holling-Tanner response function. This modification aims to capture more accurate predator-prey interactions.
\par
Estimating the parameters involved in dynamical systems is a challenging optimization problem that needs special consideration due to certain unknown natural parameters. Hence, the parameter estimation problem remains challenging despite the variety of methods available to estimate them. Artificial neural networks (ANN) \cite{ann} can be used to find the unknown parameters of the dynamical systems. "ANNs are information-processing systems with characteristics similar to biological neural networks" \cite{odeNNbook}. A neural network consists of several layers, the first layer is called the input layer, the last one is the output layer, and the layers between are called hidden layers. Every layer has several neurons each of which consists of activation function, weight value, and bias value \cite{neuron}. A deep neural network (DNN) is an ANN with more than one hidden layer. DNNs are useful for complex optimization problems due to the large number of trainable parameters. Additionally, Raissi et al. \cite{PINNpde} proposed a method called physics-informed neural network (PINN) for solving partial differential equations. In PINN, a loss function is formed using training data points, initial or boundary conditions, and residual of the differential equation. The problem is then converted into an optimization problem which can then be optimized using gradient-based optimizers like gradient descent and stochastic-based optimizers. In this study, a complex physics-informed deep neural network algorithm is employed to estimate the parameters of the proposed model using the natural population dataset.

\subsection{Motivation and Contribution}\label{motivation_and_contribution}
In the classical Lotka-Volterra model, only the searching time taken by the predators is taken into consideration but in addition to that, prey require some handling time which naturally means time required to catch and kill the prey after the search. Also in terrestrial ecology, there are scavenger species populations that can consume the dead bodies of other animals. So, this research paper aims to create and analyze a three-population predator-prey scavenger mathematical model in which respective prey response to predation is considered as Holling's functional response of type III and scavenger predator (dead bodies) interaction is considered as Holling's functional response of type I as dead bodies don't require any handling time.
\par
\textit{The objective of this research paper} is to formulate a complex non-linear dynamical mathematical model that will help us understand the interactions between predator, prey, and scavengers in ecology. The analysis of the dynamical model thus created and some examples to validate the model are also required. By developing a predator-prey scavenger model, we aim to contribute to the advancement of ecological modeling and population dynamics by adding the scavenger population in existing models along with the parameter estimation of the proposed model using the given dataset \cite{nndata}. The complex integration of DNN, PINN, and fine-tuning optimization helps efficiently estimate the parameters. Finally, the last objective of this paper is to analyze the stability of the model using the obtained parameters.
\par
\textit{This paper is organized as follows: } Section \ref{prelims} covers the preliminaries and basic concepts required for the understanding of the subsequent sections. Section \ref{ModelFormulation} covers predator-prey scavenger model formulation and Section \ref{MathematicalAnalysis} covers the mathematical analysis of the model. Section \ref{NumericalSimulations} consists of the numerical examples and respective simulations for the same. Section \ref{param_estimation} consists, of problem formulation subsection \ref{problem_formulation}, physics-informed deep neural network model subsection \ref{PIDNN_model}, and error analysis along with results in subsection \ref{error_analysis}. Section \ref{parameter_analysis} consists of stability analysis of thus found parameters using given dataset \cite{nndata}. Finally, section \ref{conclusion} concludes the paper along with expected future work.

\section{Preliminaries}\label{prelims}
This section consists of basic terminologies and optimization algorithm concepts required for the foundation of the objective of this study which is the creation of a predator-prey scavenger model and estimating natural parameters using a physics-informed deep neural network. These terminologies help to understand the subsequent section in this research paper.

\subsection{Holling's Functional Responses}\label{HollingsFunctionalResponse}
\begin{figure}[!h]
    \centering
    \includegraphics[width=0.6\textwidth]{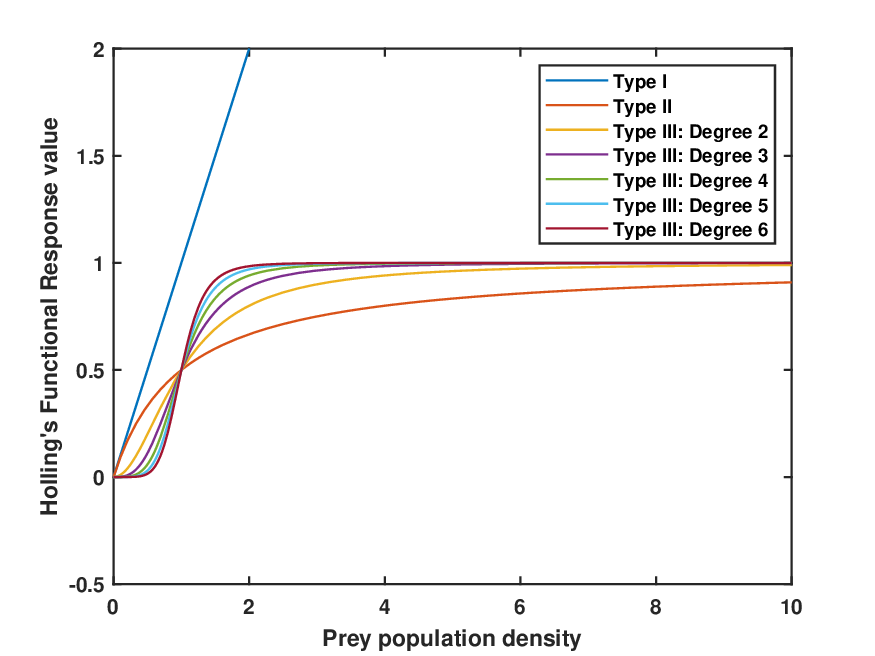}
    \caption{Holling's functional responses.}
    \label{fig:Hollings}
\end{figure}

Functional response is also very important in studying population dynamics. Functional response is majorly used in the work of Holling \cite{hollingintro, Hol1, Hol2, Hol3}. Holling’s different functional response defines the interaction between prey and predator using simple mathematical functions. Holling’s functional responses as given in Fig. \ref{fig:Hollings} are of three types,

\begin{enumerate}
    \item \textbf{Type I:} $y=ax$\\
        It is a simple function used in the classical Lotka-Volterra model and is generally used nowadays but has some natural limitations. It is linear, which shows prey don't learn anything at any prey density in response to the predator's attacks. Also, the consumption time of the predator is taken as zero, which is not natural.
    \item \textbf{Type II:} $y=\frac{ax}{1+bx}$\\
        It is a simple non-linear mathematical function that shows prey learns as prey density increases and slope decreases. Here consumption time of the predator is also taken into consideration which is defined by $b$. But in type II, learning is considered from the lowest prey density onwards. It implies prey starts to learn without any previous experience which is again not natural.
    \item \textbf{Type III:} $y=\frac{ax^n}{1+bx^n}$ where, $n\geq 2$\\
        It is the mathematical function that is the general form of type II and is significantly less taken into consideration \cite{Hol1, Hol2, Hol3} in the work of population dynamics due to its complexity. It is non-linear, which shows prey initially don’t know how to defend themselves from predators, so the slope of the graph increases up to a specific prey density. Then prey learns how to escape from the predator, and hence slope gradually decreases. Here $n$ is defined by the learning rate of the prey species.
\end{enumerate}

\subsection{Adam Algorithm}\label{adam}
Adam algorithm (Adaptive Moment Estimation) \cite{adam}, is an adaptive learning rate algorithm to improve the efficiency of the backpropagation algorithm in deep neural networks. It changes the weights and biases of the neural net based on the history of change of gradient and this helps in the prevention of avalanche effect in parameter updation. Steps involved in adam are given in algorithm \ref{algo:adam}.

\begin{algorithm}
\DontPrintSemicolon
\BlankLine
\KwIn{$\beta_1, \beta_2, g, \theta, num\_steps$}
\KwOut{$\theta$}
\BlankLine
$m_0 = 0$\;
$v_0 = 0$\;
\For{$t \gets 1$ \textbf{to} num$\_$steps}{
    $m_t = \beta_1 m_{t-1} + \left(1-\beta_1\right)g_t$\;
    $v_t = \beta_2 v_{t-1} + \left(1-\beta_2\right)g_t^2$\;
    $m_t^* = \frac{m_t}{1 - \beta_1^t}$\;
    $v_t^* = \frac{v_t}{1 - \beta_2^t}$\;
    $\theta = \theta - \frac{\alpha m_t^*}{\sqrt{v_t^* + \epsilon}}$
}
\Return{$\theta$}
\caption{{\sc Adam Algorithm}}
\label{algo:adam}
\end{algorithm}

In algorithm \ref{algo:adam}, $\alpha$ represents the step size, $m_t$ is momentum at step $t$, $g_t$ is gradient at step $t$, $\beta_1$ is decay rate of momentum, $\beta_2$ is decay rate of squared gradient, $\epsilon$ is small parameter to prevent division from zero, and $\theta$ is parameter vector.

\subsection{Broyden-Fletcher-Goldfarb-Shanno Optimization Algorithm}\label{BFGS}
The Broyden-Fletcher-Goldfarb-Shanno (BFGS) \cite{bfgs-b, bfgs-f, bfgs-g, bfgs-s} optimization algorithm is an iterative optimization method used to find the minimum of a differentiable. For easy understanding of the algorithm let's assume, We want to minimize a differentiable, scalar-valued function $f:\mathcal{R}^n\to \mathcal{R}$ with respect to input vector $x$. Initialize $x_0$ (initial guess) and $B_0$ (initial inverse Hessian matrix) and iteratively update $x$ and $B$ at each iteration $k$ according to algorithm \ref{algo:bfgs}.

\begin{algorithm}
\DontPrintSemicolon
\BlankLine
\KwIn{$x_0, B_0, f, max\_iterations$}
\KwOut{$x$}
\BlankLine
$x_1 = x_0$\;
$B_1 = B_0$\;
\For{$k \gets 1$ \textbf{to} max$\_$iterations}{
    $p_k = -B_k \nabla f\left(x_k\right)$\;
    Find $\alpha_k$ s.t. $\min f\left(x_k + \alpha_k p_k\right)$\;
    $x_{k+1} = x_k + \alpha_k p_k$\;
    $\nabla f_{k+1} = \nabla F\left(x_{k+1}\right)$\;
    $\Delta x_k = x_{k+1} - x_k$\;
    $\Delta \nabla f_k = \nabla f_{k+1} - \nabla f_k$\;
    $B_{k+1} = B_k + \frac{\Delta \nabla f_k\Delta \nabla f_k^T}{\Delta \nabla f_k^T \Delta x_k} - \frac{B_k \Delta x_k \Delta x_k^T B_k^T}{\Delta x_k^T B_k \Delta x_k}$
}
\Return{$x_k$}
\caption{{\sc BFGS Optimization Algorithm}}
\label{algo:bfgs}
\end{algorithm}

The algorithm \ref{algo:bfgs} converges to a point $x^*$ that represents the local minimum of the function $f(x)$. The BFGS optimization algorithm is widely used due to its efficiency and effectiveness in finding local minima \cite{bfgs-book}.

\section{Model Formulation}\label{ModelFormulation}
Consider three species populations, namely predator, prey, and scavengers. For modeling this system, some assumptions are required,
\begin{enumerate}
    \item The prey population (denoted by $x(t)$) can only show logistic growth and can be eaten by predators and scavengers. The only reasons for the death of the prey population are considered competition due to scarcity of resources and predation by predators and scavengers.
    \item Predator populations (denoted by $y(t)$) are considered apex predators, that is predators can eat both the prey population and scavenger population.
    \item Scavenger populations (denoted by $z(t)$) can eat the prey populations and the dead bodies of the predator populations.
\end{enumerate}
Here, $d = m_\alpha a, g = m_\beta b$ and $f = m_\gamma i$ where, $m_\alpha, m_\beta$ and $m_\gamma$ are conversion factors. The conversion factor corresponds to an increase in the predator population due to predation of single prey. 
\begin{table}[!h]
    \centering
    \begin{tabular}{|c|c|c|} \hline 
         \textbf{S. No.}&  \textbf{Parameter}& \textbf{Biological definition}\\ \hline 
         01.& $r$ &logistic growth rate of the prey population\\ \hline 
         02.& $k$ &carrying capacity of the prey population\\ \hline 
         03.& $a$ &discovery rate of prey with respect to predator\\ \hline 
         04.& $a_0$ &handling time for prey with respect to predator\\ \hline 
         05.& $b$ &discovery rate of prey with respect to scavenger\\ \hline 
         06.& $b_0$ &handling time for prey with respect to scavenger\\ \hline 
         07.& $i$ &discovery rate of scavenger with respect to predator\\ \hline 
         08.& $i_0$ &handling time for scavenger with respect to predator\\ \hline 
         09.& $d$ &rate of change of predator in presence of prey\\ \hline
         10.& $f$ &rate of change of predator in presence of scavenger\\ \hline
         11.& $g$ &rate of change of scavenger in presence of prey\\ \hline
         12.& $h$ &scavenge factor of scavengers\\ \hline
         13.& $e$ &natural death rate of predator\\ \hline
         14.& $j$ &natural death rate of scavenger\\ \hline
    \end{tabular}
    \caption{Definitions of parameters used}
    \label{tab:parameters}
\end{table}
\par
A Mathematical model to take a real situation, all the parameters must have positive values. Biological definitions of the parameters used in the model formulation are given in Table \ref{tab:parameters}. Discovery rate of prey with respect to the predator can be consider as the success rate of the predator in capturing the prey. Handling time for prey with respect to predator is the time required by the predator to catch and kill the prey after searching for it.
\par
\begin{figure}[!h]
    \centering
    \includegraphics[scale=0.6]{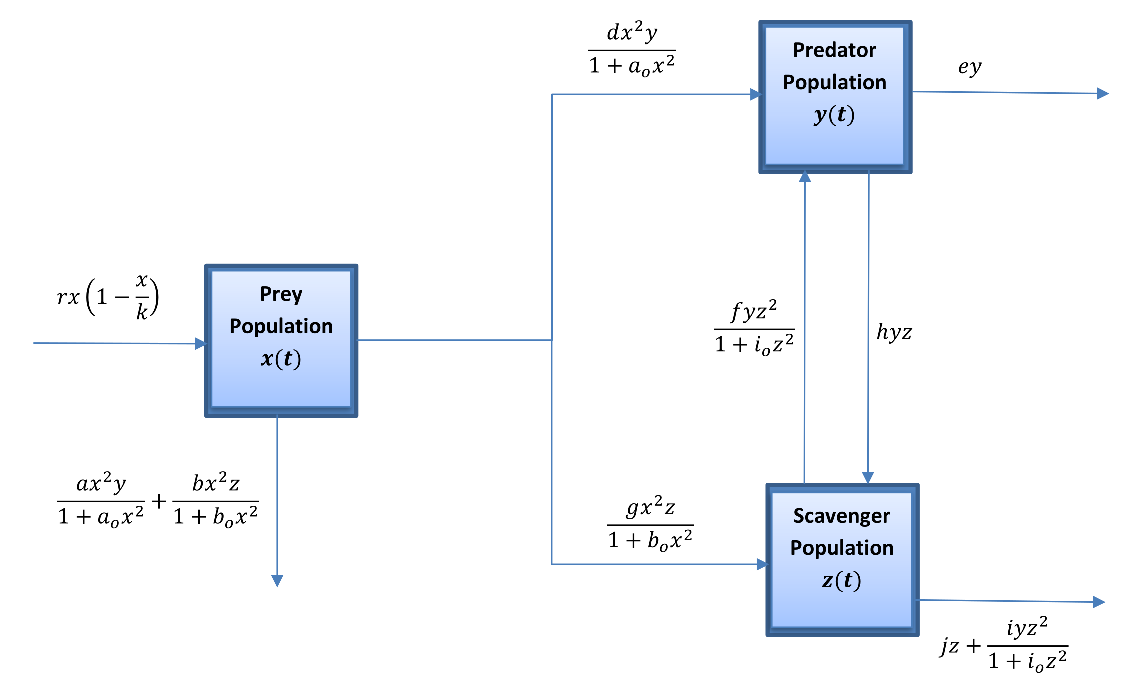}
    \caption{Predator prey scavenger interaction flowchart}
    \label{fig:interaction}
\end{figure}

Now, the proposed mathematical predator prey scavenger model (Fig. \ref{fig:interaction}) is,
\begin{equation}\label{eq:propose}
\begin{split}    
   \frac{dx}{dt} &= rx\left(1 - \frac{x}{k}\right) - \frac{ax^2 y}{1 + a_0 x^2} - \frac{bx^2 z}{1 + b_0 x^2}\\
 \frac{dy}{dt} &= \frac{dx^2 y}{1 + a_0 x^2} + \frac{fz^2 y}{1 + i_0 z^2} - ey\\
 \frac{dz}{dt} &= \frac{gx^2 z}{1 + b_0 x^2} + hyz - \frac{iyz^2}{1 + i_0 z^2} - jz 
\end{split}
\end{equation}
with initial conditions as $x(0)>0,$ $y(0)>0$ and $z(0)>0.$
\par
In the proposed mathematical model, prey-predator, prey-scavenger, and predator-scavenger interactions (predation interaction only) are taken as Holling's functional response of type III as prey and scavengers require some handling time with respect to their respective predators. Also, they start to learn to defend themselves as the respective prey density increases. But scavenger-predators (dead) interaction is taken as Holling's functional response of type I as dead bodies of predator species don't require any handling time which is obvious.
\par
In the proposed model, Holling's functional response of type III with $n=2$ is considered as it provides the learning curve of the respective prey species more realistic because it gradually increases the slope and then gradually decreases it also. Due to gradual learning of brain with experience, $n=2$ curve seems to be realistic. A higher degree $n>2$ is not reasonable for the natural systems as the learning curve did not seems realistic because $n>2$ shows that prey species learns at a much faster rate which is not realistic.

\section{Mathematical Analysis of Model}\label{MathematicalAnalysis}
In this section, analysis of different systems is done, which are extracted from the predator prey scavenger model \eqref{eq:propose}. There are many systems which can be derived, given in Table \ref{tab:DiffSystems}. Out of systems given in Table \ref{tab:DiffSystems}, prey only, predator only and scavenger only are trivial systems.
\begin{table}[!h]
    \centering
    \begin{tabular}{|c|c|c|c|} \hline 
          \textbf{S. No.}&\textbf{Model}& \textbf{Consideration} &\textbf{Remarks} \\ \hline 
          01.&Prey only model& $y(t)=0$ and $z(t)=0$ &exist in nature\\ \hline 
          02.&Predator only model& $x(t)=0$ and $z(t)=0$ &can't exist in nature\\ \hline 
          03.&Scavenger only model& $x(t)=0$ and $y(t)=0$ &can't exist in nature\\ \hline 
          04.&Predator scavenger model& $x(t)=0$ &exist in nature\\ \hline 
          05.&Scavenger prey model& $y(t)=0$ &exist in nature\\ \hline 
          06.&Predator prey model& $z(t)=0$ &exist in nature\\ \hline 
          07.&Predator prey scavenger model&  -&exist in nature\\ \hline
    \end{tabular}
    \caption{Different systems derived from predator prey scavenger model}
    \label{tab:DiffSystems}
\end{table}

\subsection{Predator Scavenger Model}\label{PredScav}
Scavengers can only feed on some proportion of dead bodies of the predators, and predators can feed on scavengers only. As population of scavengers decreases due to predation, predator population also decreases, resulting in a much faster decrease in scavenger population. Hence, predator scavengers can't exist in nature without prey which can be seen from the analysis of the predator scavengers model,
\begin{equation}\label{eq:PredScav}
\begin{split}
    \frac{dy}{dt} &= \frac{fz^2 y}{1 + i_0 z^2} - ey\\
    \frac{dz}{dt} &= hyz - \frac{iyz^2}{1 + i_0 z^2} - jz 
\end{split}
\end{equation}
The steady state points are, $(0,0)$ and $\left( z_0 , \frac{fj z_0}{hfz_0-ie}\right)$ where, $z_0 = \left[ \frac{e}{f-i_0 e} \right]^{1/2}$.
\begin{enumerate}
    \item \textit{\textbf{Steady state}} $(0,0)$,\\
        The jacobian matrix for this case is,
        $$J=\begin{bmatrix}-e&0\\ 0&-j\end{bmatrix} $$
        implies eigenvalues are $-e$ and $-j$, which means steady state $(0,0)$ is stable as both eigenvalues are negative.
    \item \textit{\textbf{Steady state}} $\left( z_0 , \frac{fjz_0}{hfz_0 - ie}\right)$,\\
        The existence conditions for this steady state point are $f-i_0e > 0$ and $hfz_0 - ie > 0$. The jacobian matrix for predator scavenger model for this steady point case is,
        $$J_{yz} = 
        \begin{bmatrix}
        0&\frac{2je(f-i_0e)}{hfz_0 - ie} \vspace{0.2in}\\ 
        hz_{0}-\frac{ie}{f} \hspace{0.2in}& \frac{jei(2i_0e-f)}{f(hfz_0 - ie)}
        \end{bmatrix}$$
        Hence, the characteristic equation is given as,
        $$m(\lambda) = \lambda^{2} -\lambda \left[\frac{jei(2i_0e-f)}{f(hfz_0-ie)} \right]  -\frac{2je(f-i_0e)}{f} = 0 $$
        From the above equation, the eigenvalues ($\lambda_1$ and $\lambda_2$) gives, 
        $$\lambda_1\lambda_2 = -\frac{2je(f-i_0e)}{f}<0$$
        According to Routh Hurwitz criterion \cite{routh-hur}, steady state $\left( z_0 , \frac{fjz_0}{hfz_0-ie}\right)$ is unstable.
\end{enumerate}

\subsection{Predator Prey Model}\label{PredPrey}
Predator prey without scavengers can exist in nature which can be seen from the analysis of the predator prey model,
\begin{equation}\label{eq:PredPrey}
\begin{split}
    \frac{dx}{dt} &= rx\left(1 - \frac{x}{k}\right) - \frac{ax^2 y}{1 + a_0 x^2}\\
    \frac{dy}{dt} &= \frac{dx^2 y}{1 + a_0 x^2} - ey 
\end{split}
\end{equation}
The steady state points are $(0,0)$, $(k,0)$ and $\left( x_0 , \frac{drx_0}{aek} (k-x_0)\right)$ where, $x_0 = \left[ \frac{e}{d-a_0 e} \right]^{1/2}$.
\begin{enumerate}
    \item \textbf{\textit{Steady state}} $(0,0)$,\\
        The jacobian matrix for this case is,
        $$J=\begin{bmatrix} r&0\\ 0&-e\end{bmatrix} $$
        that is eigenvalues are $r$ and $-e$, which means $(0,0)$ is an unstable steady state point as one of the eigenvalue is positive.
    \item \textbf{\textit{Steady state}} $(k,0)$,\\
        The jacobian matrix for this case is,
        $$ J_{x}=\begin{bmatrix}-r&\frac{-a k^2}{1+a_0 k^2} \vspace{0.1in}\\ 0&\hspace{0.1in}\frac{dk^2}{1+a_0k^2} -e\end{bmatrix} $$
        that is this steady state point is stable if $k^2 < \frac{e}{d-a_0 e}$.
    \item \textbf{\textit{Steady state}} $\left( x_0 , \frac{drx_0}{aek} (k-x_0)\right)$,\\
        The existence condition for this steady state point is $d-a_0e > 0$. The jacobian matrix for this case is,
        $$J_{xy} = 
        \begin{bmatrix}
        -r+\frac{2ra_0 e (k-x_0)}{d k} &\frac{-ae}{d} \vspace{0.1in}\\ 
        \frac{2r(d-a_0 e)(k-x_0)}{a k} \hspace{0.2in}&0
        \end{bmatrix}$$
        The characteristic equation is given as,
        $$m(\lambda) := \lambda^{2} -\lambda \left[ -r+\frac{2ra_0 e (k-x_0)}{d k} \right]  +\frac{2re(d-a_0 e)(k-x_0)}{d k} = 0\text{.}$$
        Since, population capacity $>$ population at any instantaneous point $(k>x_0)$. It is clear from above that $m(0) > 0$ and using Routh Hurwitz criterion \cite{routh-hur}, $\lambda_1\lambda_2 >0$. Therefore, either both roots are positive or both negative. For both roots ($\lambda_1 \& \lambda_2$) to be negative, $\lambda_1\lambda_2 >0$ and $\lambda_1 + \lambda_2 < 0$. Hence,
        $$\lambda_1 + \lambda_2 = \frac{2ra_0 e (k-x_0)}{d k} - r < 0\text{.}$$
        So, the steady state $\left( x_0 , \frac{drx_0}{aek} (k-x_0)\right)$ is stable if $\frac{2a_0 e (k-x_0)}{d k} < 1$.
\end{enumerate}

\subsection{Scavenger Prey Model}\label{ScavPrey}
Scavenger prey without predators can exist in nature (because scavenger prey interaction is similar to predator prey interaction) which can be seen from the analysis of the scavenger prey model,
\begin{equation}\label{eq:ScavPrey}
\begin{split}
    \frac{dx}{dt} &= rx\left(1 - \frac{x}{k}\right) - \frac{bx^2 z}{1 + b_0 x^2}\\
    \frac{dz}{dt} &= \frac{gx^2 z}{1 + b_0 x^2} - jz 
\end{split}
\end{equation}
The steady state points are $(0,0)$, $(k,0)$ and $\left( x_0 , \frac{grx_0}{bjk} (k-x_0)\right)$ where, $x_0 = \left[ \frac{j}{g-b_0 j} \right]^{1/2}$.
\begin{enumerate}
    \item \textbf{\textit{Steady state}} $(0,0)$,\\
        The jacobian matrix for this case is,
        $$J=\begin{bmatrix} r&0\\ 0&-j\end{bmatrix} $$
        that is eigenvalues are $r$ and $-j$ which means $(0,0)$ is an unstable steady state point as one of the eigenvalue is positive.
    \item \textbf{\textit{Steady state}} $(k,0)$,\\
        The jacobian matrix for this case is,
        $$ J_{x}=\begin{bmatrix}-r&\frac{-b k^2}{1+b_0 k^2} \vspace{0.1in}\\ 0& \hspace{0.2in}\frac{gk^2}{1+b_0k^2} -j\end{bmatrix} $$
        that is this steady state point is stable if $k^2 < \frac{j}{g-b_0 j}$.
    \item \textbf{\textit{Steady state}} $\left( x_0 , \frac{grx_0}{bjk} (k-x_0)\right)$,\\
        The existence condition for this steady state point is $g-b_0j > 0$. The jacobian matrix for this case is,
        $$J_{xz} = 
        \begin{bmatrix}
        -r+\frac{2rb_0 j (k-x_0)}{g k} &\frac{-bj}{g} \vspace{0.1in}\\ 
        \frac{2r(g-b_0 j)(k-x_0)}{b k} \hspace{0.1in}&0
        \end{bmatrix}$$
        The characteristic equation is given as,
        $$m(\lambda) := \lambda^{2} -\lambda \left[ -r+\frac{2rb_0 j (k-x_0)}{g k} \right]  +\frac{2rj(g-b_0 j)(k-x_0)}{g k} = 0\text{.}$$
        Since, population capacity $>$ population at any instantaneous point $(k>x_0)$. It is clear from above that $m(0) > 0$ and using Routh Hurwitz criterion \cite{routh-hur}, $\lambda_1\lambda_2 >0$. Therefore, either both roots are positive or both negative. For both roots ($\lambda_1 \& \lambda_2$) to be negative, $\lambda_1\lambda_2 >0$ and $\lambda_1 + \lambda_2 < 0$. Hence,
        $$\lambda_1 + \lambda_2 = \frac{2rb_0 j (k-x_0)}{g k} - r < 0\text{.}$$
        So, the steady state $\left( x_0 , \frac{grx_0}{bjk} (k-x_0)\right)$ is stable if $\frac{2b_0 j (k-x_0)}{g k} < 1$.
\end{enumerate}

\subsection{Predator Prey Scavenger Model}\label{PredPreyScav}
Considering all populations in the system viz. prey, predator and scavenger. Predator prey scavengers all can coexist in nature which can be seen from the analysis of the predator prey scavenger model,
\begin{equation}\label{eq:PredPreyScav}
\begin{split}
    \frac{dx}{dt} &= rx\left(1 - \frac{x}{k}\right) - \frac{ax^2 y}{1 + a_0 x^2} - \frac{bx^2 z}{1 + b_0 x^2}\\
    \frac{dy}{dt} &= \frac{dx^2 y}{1 + a_0 x^2} + \frac{fz^2 y}{1 + i_0 z^2} - ey\\
    \frac{dz}{dt} &= \frac{gx^2 z}{1 + b_0 x^2} + hyz - \frac{iyz^2}{1 + i_0 z^2} - jz 
\end{split}
\end{equation}
Let steady states for the system be denoted by $(0,0,0)$, $(E_x^{[2]},0,0)$, $(0,E_y^{[3]},E_z^{[3]})$, $(E_x^{[4]},E_y^{[4]},0)$, $(E_x^{[5]},0,E_z^{[5]})$ and $(E_x^{[6]},E_y^{[6]},E_z^{[6]})$. It is clear from Table \ref{tab:DiffSystems}, the cases of ($0,E_y,0$) and ($0,0,E_z$) are included in ($0,0,0$) case as both yields same stable steady state point $(0,0,0)$.
\begin{enumerate}
    \item \textbf{\textit{Steady state}} $(0,0,0)$,\\
        The jacobian matrix for this case is,
        $$J_0 = 
        \begin{bmatrix}
        r&0&0\\ 0&-e&0\\ 0&0&-j
        \end{bmatrix} $$
        that is eigenvalues are $r,-e$ and $-j$, which means $(0,0,0)$ is an unstable steady state point as one of the eigenvalue is positive.
    \item \textbf{\textit{Steady state}} $(E_x^{[2]},0,0) = (k,0,0)$,\\
        The jacobian matrix for this case is,
        $$J_x = 
        \begin{bmatrix}
        -r& \hspace{0.1in} \frac{-a k^2}{1+a_0 k^2}&\hspace{0.1in} \frac{-b k^2}{1+b_0 k^2} \vspace{0.1in} \\
        0\hspace{0.1in}&\frac{dk^2}{1+a_0 k^2} -e&\hspace{0.1in} 0 \vspace{0.1in}\\ 
        0\hspace{0.1in} &0 \hspace{0.1in}&\frac{g k^2}{1+b_0 k^2}-j
        \end{bmatrix} $$
        that is steady state $(E_x^{[2]},0,0)$ is stable if $k^2 < \frac{e}{d-a_0 e}$ and $k^2 < \frac{j}{g-jb_0}$.
    \item \textbf{\textit{Steady state}} $(0,E_y^{[3]},E_z^{[3]}) = \left( 0,\frac{fjz_{0}}{hfz_0 - ie} ,z_{0}\right)$ where, $z_0 = \left[ \frac{e}{f-i_{0}e} \right]^{1/2}$\\
        The existence conditions of this steady state are $f-i_0 e >0$ and $hfz_0 - ie>0$. The jacobian matrix for this case is,
        $$J_{yz} = 
            \begin{bmatrix}
            r&0&0\vspace{0.1in}\\ 
            0&0&\frac{2je(f-i_0 e)}{hfz_0 -i e} \vspace{0.1in}\\ 
            0\hspace{0.1in}&hz_0 - \frac{ie}{f} &\hspace{0.1in}\frac{jei(2i_0e-f)}{f(hfz_0 - ie)}
            \end{bmatrix} $$
        which means steady state $(0,E_y^{[3]},E_z^{[3]})$ is unstable as one of the eigenvalues is $r$ which is positive.
    \item \textbf{\textit{Steady state}} $(E_x^{[4]},E_y^{[4]},0) = \left( x_0,\frac{drx_{0}}{aek} (k-x_0),0\right)$ where, $x_0 = \left[ \frac{e}{d-a_{0}e} \right]^{1/2} $\\
        The existence condition of this steady state is $d-a_0e>0$. The jacobian matrix for this case is,
        $$J_{xy} = 
            \begin{bmatrix}
            -r+ \frac{2ra_0 e (k-x_0)}{d k} \hspace{0.1in}&\frac{-a e}{d} & \frac{-be}{d+(b_0 - a_0)e}\vspace{0.1in}\\
            \frac{2r(d-a_0e)(k-x_0)}{ak}&0&0\vspace{0.1in}\\
            0&0&\hspace{0.1in}\frac{ge}{d+(b_0-a_0)e}+\frac{hdrx_0(k-x_0)}{aek}-j
            \end{bmatrix}$$
        This steady state is stable if,
        \begin{enumerate}
            \item[\textbf{1.}] $\frac{2a_0 e (k-x_0)}{d k} < 1$ (from Section \ref{PredPrey})
            \item[\textbf{2.}] $hdrx_0\left(1-\frac{x_0}{k}\right)[d+(b_0 - a_0)e] + ae[ge-j(d+(b_0 - a_0)e)] <0$
        \end{enumerate}
        As condition 2. corresponds to one of the eigenvalues, hence must be negative for the system to be stable.
    \item \textbf{\textit{Steady state}} $(E_x^{[5]},0,E_z^{[5]}) = \left( x_0,0, \frac{grx_{0}}{bjk} (k-x_0)\right)$ where, $x_0 = \left[ \frac{j}{g-jb_{0}} \right]^{1/2} $\\
        The existence condition of this steady state is $g-b_0j>0$. The jacobian matrix for this case is,
        $$J_{xz} =
            \begin{bmatrix}
            -r + \frac{2rb_0 j (k-x_0)}{g k} \hspace{0.1in}&\frac{-a j}{g+j(a_0-b_0)} &\hspace{0.1in} \frac{-bj}{g} \vspace{0.1in}\\
            0 \hspace{0.1in}& \frac{dj}{g+j(a_0 - b_0)} +\frac{fz_0^2}{1+i_0 z_0^2} -e \hspace{0.1in}& 0 \vspace{0.1in}\\
            \frac{2r(g-jb_0)(k-x_0)}{bk} \hspace{0.1in}& \frac{hgrx_0(k-x_0)}{bjk} - \frac{iz_0^2}{1+i_0z_0^2} \hspace{0.1in}& 0
            \end{bmatrix}$$
        where, $z_0 = \frac{grx_0}{bjk}(k-x_0)$. This steady state is stable if,
        \begin{enumerate}
            \item[\textbf{1.}] $\frac{2b_0 j (k-x_0)}{g k} < 1$ (from Section \ref{ScavPrey})
            \item[\textbf{2.}] $\frac{fz_0^{2}}{1+i_0 z_0^{2}}< \frac{e(g+j(a_0 -b_0))-dj}{g+j(a_0 -b_0)}~~$($t$(say))$ \Rightarrow ~~z_0^{2} < \frac{t}{f-i_0 t}$
        \end{enumerate}
        As condition 2. corresponds to one of the eigenvalues, hence must be negative for the system to be stable.
    \item \textbf{\textit{Steady state}} $(E_x^{[6]},E_y^{[6]},E_z^{[6]}) = (x^*,y^*,z^*)$\\
        where $x^*$ is the only positive real root of equation (see Appendix for values of $p_{i}$),
        \begin{equation}\label{eq:appendix}
           \sum^{12}_{i=0} p_{i}X^{i}=0 
        \end{equation}
        The values of $y^*$ and $z^*$ are given as,
        \begin{equation*}
        \begin{split}
            {z^*}^2 &= \frac{e+(a_0 e - d)x^{*^2}}{f(1+a_0 x^{*^2})-i_0 \left(e+(a_0 e - d)x^{*^2}\right)}\\
            y^* &= \left( j- \frac{gx^{*^2}}{1+b_0 x^{*^2}}\right)\left(\frac{1+i_0 z^{*^2}}{h+i_0 z^{*^2} - iz^*}\right)
        \end{split}
        \end{equation*}
        The conditions for existence of this steady state are,
        \begin{enumerate}
            \item[\textbf{1.}] $x^*$ must be only positive real root of equation $\sum_{i=0}^{12} p_{i}X^{i}=0$
            \item[\textbf{2.}] $z^* < min\left\{\frac{r(1+b_0x^{*^2})(k-x^*)}{bx^*}, \frac{ie}{hf}\right\}$
        \end{enumerate}
        The jacobian matrix for this case is,
        $$J_{xyz} = 
        \begin{bmatrix}
            r(1-\frac{2x^*}{k}) - \frac{2ax^*y^*}{(1+a_0 x^{*^2})^2} - \frac{2bx^*z^*}{(1+b_0 x^{*^2})^2} \hspace{0.091in}& \frac{-ax^{*^2}}{1+a_0 x^{*^2}} & \frac{-bx^{*^2}}{1+b_0 x^{*^2}}\vspace{0.1in} \\
            \frac{2dx^*y^*}{(1+a_0 x^{*^2})^2} & 0 & \frac{2fz^*y^*}{(1+i_0 z^{*^2})^2}\vspace{0.1in}\\
            \frac{2gx^*z^*}{(1+b_0 x^{*^2})^2} & hz^* - \frac{iz^{*^2}}{1+i_0 z^{*^2}}\hspace{0.09in} & \frac{iy^*z^*}{1+i_0 z^{*^2}} - \frac{2iy^*z^*}{(1+i_0 z^{*^2})^2}
        \end{bmatrix}$$
        Considering the above matrix as,
        $$J = 
        \begin{bmatrix}
        a_{11} & a_{12} & a_{13}\\
        a_{21} & a_{22} & a_{23}\\
        a_{31} & a_{32} & a_{33}
        \end{bmatrix}
        $$Let's assume some variables as,
        \begin{enumerate}
            \item[\textbf{1.}] $m_1 = -trace(J)$
            \item[\textbf{2.}] $m_2 = M_{11} + M_{22} + M_{33} \quad$ where, $M_{ii}$ is the minor of $a_{ii}^{th}$ element, for $i=1,2,3$.
            \item[\textbf{3.}] $m_3 = -determinant(J)$
        \end{enumerate}
        These variables are such that the characteristic equation of matrix J is, $\lambda^3 + m_1\lambda^2 + m_2\lambda + m_3 = 0$. According to Routh Hurwitz criterion \cite{routh-hur}, steady state $(E_x^{[6]},E_y^{[6]},E_z^{[6]})$ stable if,
        \begin{enumerate}
            \item[\textbf{1.}] $m_i > 0$ where, $i=1,2,3$
            \item[\textbf{2.}] $m_1m_2 - m_3 > 0$
        \end{enumerate}
\end{enumerate}

\section{Numerical Simulations}\label{NumericalSimulations}
\begin{example}\label{ex1}
    Consider predator scavenger system with parameter values as  $ f=0.5, i_0 = 0.25, e=0.5, h=0.5, i=0.5$ and $j=1.5$, with initial conditions $y(0) = 4$ and $z(0) = 6$.
    \begin{equation*}
        \begin{split}
            \frac{dy}{dt} &= \frac{0.5z^2 y}{1 + 0.25 z^2} - 0.5y\\
            \frac{dz}{dt} &= 0.5yz - \frac{0.5yz^2}{1 + 0.25 z^2} - 1.5z
        \end{split}
    \end{equation*}
\end{example}

\begin{figure}[!h]
    \centering
    \begin{subfigure}[b]{7.5cm}
         \centering
         \includegraphics[width=\textwidth]{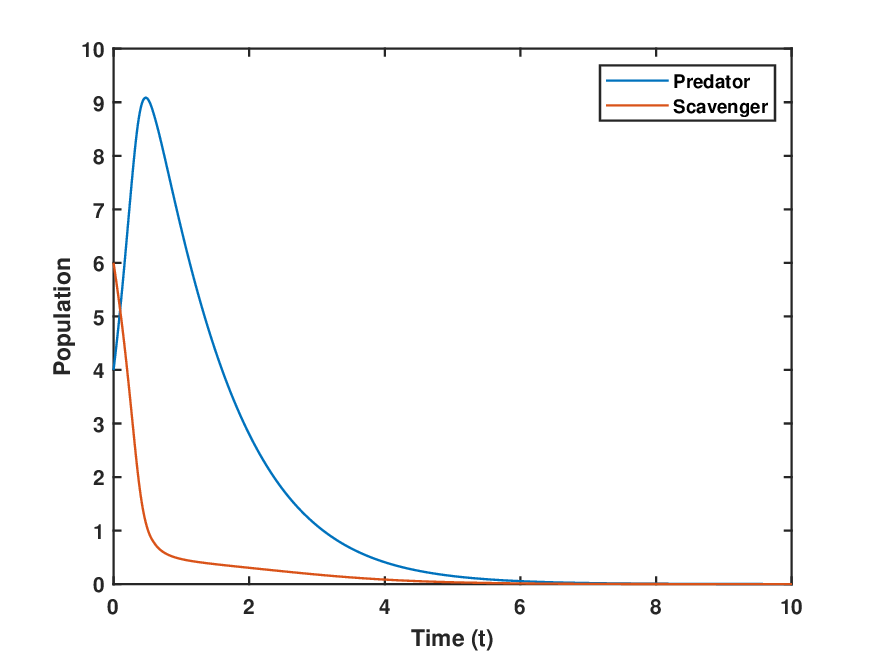}
         \caption{Dynamical behaviour.}
    \end{subfigure}
    \hfill
    \begin{subfigure}[b]{7.5cm}
         \centering
         \includegraphics[width=\textwidth]{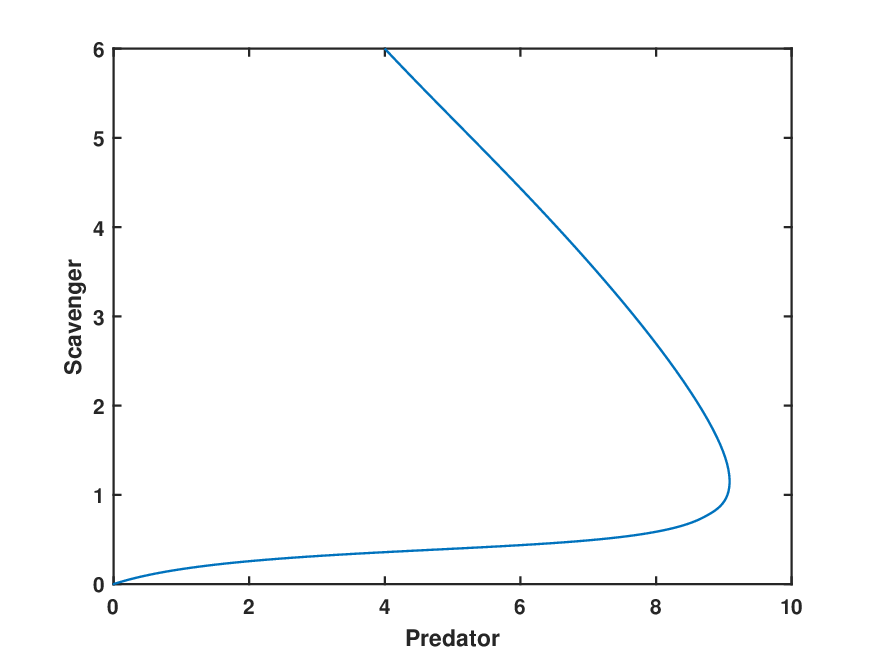}
         \caption{Phase portrait.}
    \end{subfigure}
        \caption{Predator scavenger model for Example \ref{ex1}.}
        \label{fig:ex1}
\end{figure}

For the predator scavenger system, it is clear from Fig. \ref{fig:ex1} that both populations stabilize at point $(0,0)$ in the system and converge to zero. Here, 
$$f-i_0 e = 0.375 > 0$$ 
$$~hf\left( \frac{e}{f-i_0 e} \right)^{1/2} - ie = \frac{2-\sqrt{3}}{4\sqrt{3}} = 0.038675 > 0$$ 
that means the existence condition is satisfied for the equilibrium points. Here,
$$z_0 = \left[ \frac{e}{f-i_0 e} \right]^{1/2} = \frac{2}{\sqrt{3}}$$ 
and hence, equilibrium points are $(0,0)$ and $(1.1547,22.3923)$,
$$\left( z_0 , \frac{fj z_0}{hfz_0-ie}\right) = \left(\frac{2}{\sqrt{3}} , 6(2 + \sqrt{3})\right) = (1.1547,22.3923)$$
According to subsection \ref{PredScav}, the equilibrium point $(0,0)$ is stable, and the equilibrium point $(1.1547,22.3923)$ is unstable. It is clear from the Fig. \ref{fig:ex1} that both the scavenger population as well as predator population converge to zero. The stability of $(0,0)$ can be explained easily as the only source of food for the predators is scavengers, and the only source of food for the scavengers population are the dead bodies of the predators. Initially, there are more scavengers but their natural death rate is greater, so their population decreases and the population of predators increases. As the population of scavengers decreases, food for the predators starts to decline as well. Hence, predator population also starts falling. As a result, both populations stabilizes at $(0,0)$ or, in other words both populations extinct.

\begin{example}\label{ex2}
    Consider predator prey system with parameter values as $r=1$, $k = 2$, $a=1$, $a_0=0.25$, $d=1$ and $e=1$ with initial conditions $x(0) = 2$ and $y(0) = 4$.
    \begin{equation*}
        \begin{split}
            \frac{dx}{dt} &= x\left(1 - \frac{x}{2}\right) - \frac{x^2 y}{1 + 0.25x^2}\\
            \frac{dy}{dt} &= \frac{x^2 y}{1 + 0.25x^2} - y
        \end{split}
    \end{equation*}
\end{example}

\begin{figure}[!h]
    \centering
    \begin{subfigure}[b]{7.5cm}
         \centering
         \includegraphics[width=\textwidth]{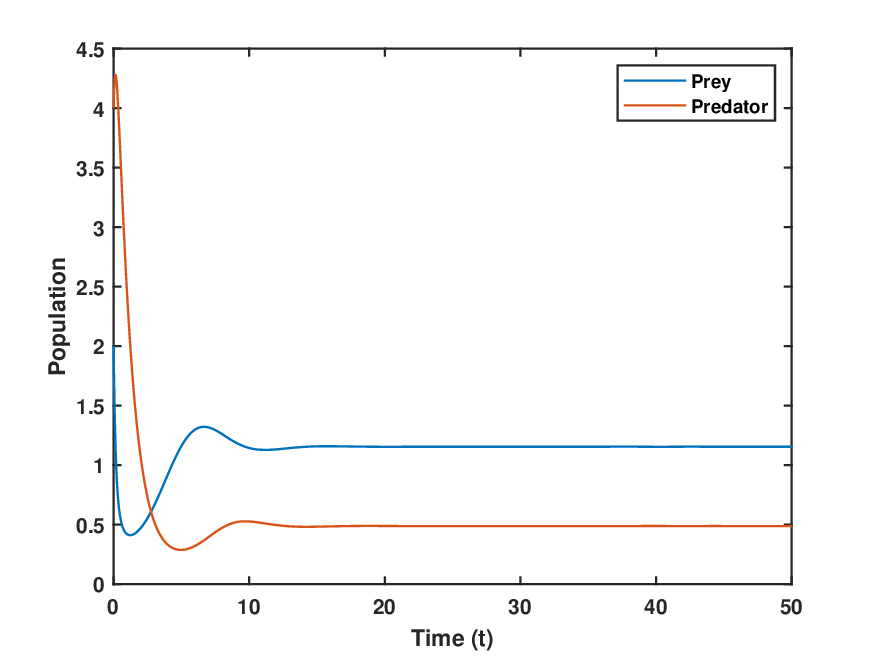}
         \caption{Dynamical behaviour.}
    \end{subfigure}
    \hfill
    \begin{subfigure}[b]{7.5cm}
         \centering
         \includegraphics[width=\textwidth]{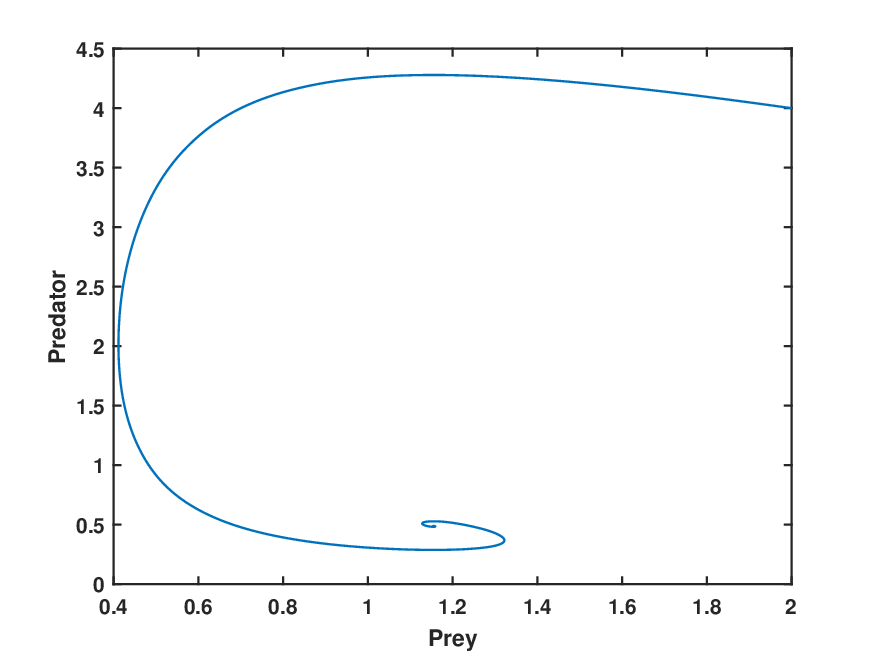}
         \caption{Phase portrait.}
    \end{subfigure}
        \caption{Predator prey model for Example \ref{ex2}.}
        \label{fig:ex2}
\end{figure}

For predator prey system, it is clear from Fig. \ref{fig:ex2} that both of the populations stabilize at some specific value of the population which is neither $(0,0)$ nor $(k,0) = (2,0)$. Here, $d-a_0 e = 0.75 > 0$ that means existence condition is satisfied for the equilibrium points. Here, $x_0 = \left[ \frac{e}{d-ea_0} \right]^{1/2} = \frac{2}{\sqrt{3}}$ and hence, equilibrium points are $(0,0)$, $(2,0)$ and $(1.1547,0.4880)$,
$$\left( x_0,\frac{drx_0}{aek} (k-x_0)\right)  = (1.1547,0.4880)$$
As discussed in the predator prey model subsection \ref{PredPrey}, the equilibrium point $(0,0)$ is unstable. For the equilibrium point $(2,0)$, $k^2 = 4$ and $\frac{e}{d-ea_0} = \frac{4}{3}$. As $k^2 > \frac{e}{d-ea_0}$, so equilibrium point $(2,0)$ is not stable because if this point gets stabilize then it means that only prey exists in nature and predators extinct which is not possible. For the equilibrium point $(1.1547,0.4880)$, $\frac{2a_0 e (k-x_0)}{d k} = 0.211325 < 1$ which satisfies the stability condition for the equilibrium point $(1.1547,0.4880)$. Hence, the system stabilizes and both predator and prey populations coexist in this system.\\
\textbf{\textit{Remarks:}} Similarly, the existence of solution and the stability condition for the scavenger prey model can be determined. In scavenger prey model, scavengers will show the same behaviour as the predator for the prey.

\begin{example}\label{ex3}
    Consider predator prey scavenger system with parameter values such as $r = 0.5$, $k = 100$, $a = 0.5$, $a_0 = 0.25$, $b = 0.5$, $b_0 = 0.25$, $d = 0.5$, $e = 1$, $f = 0.1$, $g = 0.5$, $h = 0.1$, $i = 0.1$, $i_0 = 0.25$ and $j=1$ with initial conditions $x(0) = 4$, $y(0) = 3$ and $z(0) = 2$.
    \begin{equation*}
        \begin{split}
            \frac{dx}{dt} &= 0.5x\left(1 - \frac{x}{100}\right) - \frac{0.5x^2 y}{1 + 0.25 x^2} - \frac{0.5x^2 z}{1 + 0.25 x^2}\\
            \frac{dy}{dt} &= \frac{0.5x^2 y}{1 + 0.25 x^2} + \frac{0.1z^2 y}{1 + 0.25 z^2} - y\\
            \frac{dz}{dt} &= \frac{0.5x^2 z}{1 + 0.25 x^2} + 0.1yz - \frac{0.1yz^2}{1 + 0.25 z^2} - z
        \end{split}
    \end{equation*}
\end{example}

\begin{figure}[!h]
    \centering
    \begin{subfigure}[b]{7.5cm}
         \centering
         \includegraphics[width=\textwidth]{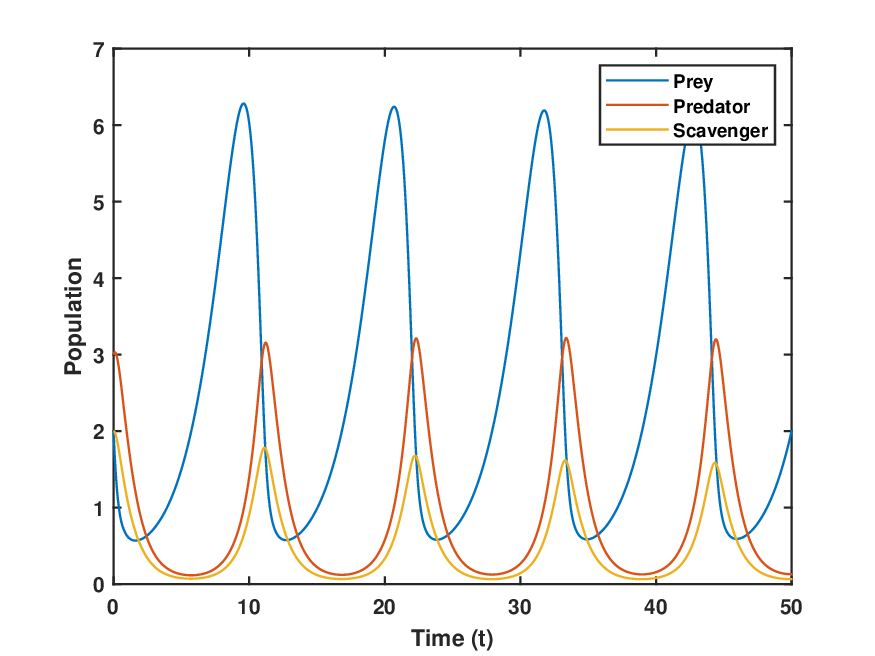}
         \caption{Dynamical behaviour.}
    \end{subfigure}
    \hfill
    \begin{subfigure}[b]{7.5cm}
         \centering
         \includegraphics[width=\textwidth]{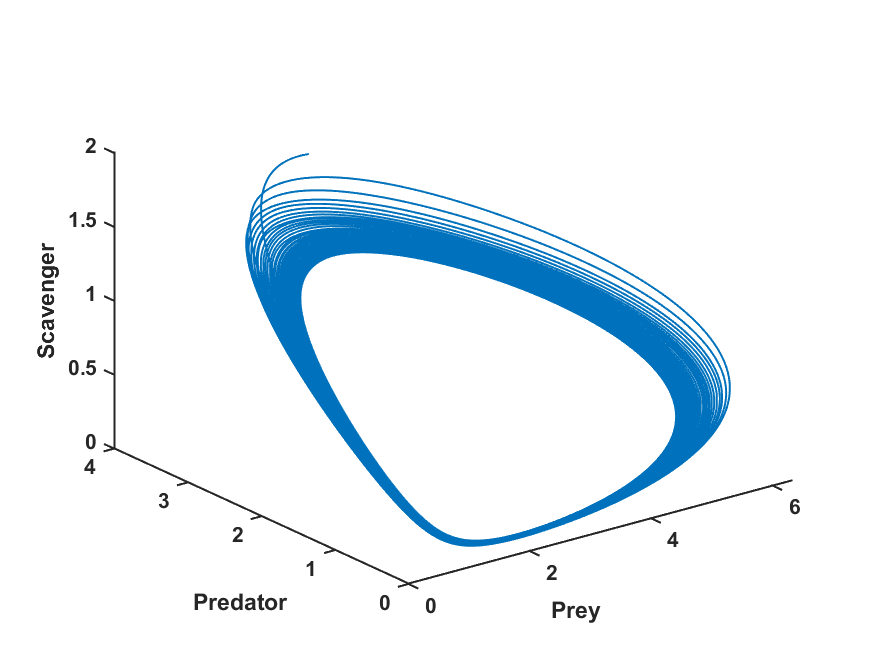}
         \caption{Phase portrait.}
    \end{subfigure}
        \caption{Predator prey scavenger model for Example \ref{ex3}.}
        \label{fig:ex3}
\end{figure}

For predator prey scavenger system, it is clear from Fig. \ref{fig:ex3} that all of the population don't stabilize at any value of the population.
\\
\textbf{\textit{Equilibrium point}} ($0,0,0$): It is clear that $(0,0,0)$ is an unstable equilibrium point because if this point becomes stable, then the system collapse. As described in subsection \ref{PredPreyScav}, ($0,E_y,0$) and ($0,0,E_z$) cases are also included in this ($0,0,0$) case as both yields equilibrium point $(0,0,0)$ and both are not stable as well.
\\
\textbf{\textit{Equilibrium point}} ($E_x^{[2]},0,0$): For equilibrium point $(100,0,0)$, 
$$\frac{e}{d-ea_0} = 4 < 100^2 = k^2 \quad\text{and}\quad \frac{j}{g-jb_0} = 4 < 100^2 = k^2$$ 
Hence, the equilibrium point $(100,0,0)$ is not stable.
\\
\textbf{\textit{Equilibrium point}} ($0,E_y^{[3]},E_z^{[3]}$): As, $f-i_0 e = -0.15 < 0$ didn't satisfy the existence condition for the equilibrium point $\left(0,\frac{fj z_0}{hfz_0-ie},z_0\right)$ because equilibrium point turns out to have a negative or, complex numbers which has no significance in nature. Hence, equilibrium point ($0,E_y^{[3]},E_z^{[3]}$) does not exist in the system.
\\
\textbf{\textit{Equilibrium point}} ($E_x^{[4]},E_y^{[4]},0$): As $d-a_0 e = 0.25 > 0$ implies the existence condition for the equilibrium point, 
$$\left( \left[ \frac{e}{d-a_0 e} \right]^{1/2} , \frac{drx_0(k-x_0)}{aek}, 0\right) = (2,0.98,0)$$
For this equilibrium point, 
$$\frac{2a_0 e (k-x_0)}{d k} = 0.98 < 1$$ 
$$hrdx_0\left(1-\frac{x_0}{k}\right)(d+(b_0-a_0)e) + ae(ge - j(d+(b_0-a_0)e)) = 0.0254 > 0$$ 
Hence, the stability conditions are not met which means predator and prey may coexist in an unstable system.
\\
\textbf{\textit{Equilibrium point}} ($E_x^{[5]},0,E_z^{[5]}$): As $g-b_0 j = 0.25 > 0$ implies the existence condition for the equilibrium point, 
$$\left( \left[ \frac{j}{g-b_0 j} \right]^{1/2} , 0 , \frac{grx_0(k-x_0)}{bjk}\right) = (2,0,0.98)$$ 
For this equilibrium point, 
$$\frac{2b_0j(k-x_0)}{gk} = 0.98 < 1$$ 
$$t = \frac{e(g+j(a_0-b_0))-dj}{g+j(a_0-b_0)} = 0$$ 
and hence $z_0 ^2 = 0.9604 > \frac{t}{f-i_0 t} = 0$. Therefore, the stability condition doesn't satisfy which means scavenger and prey may coexist in an unstable system.
\\
\textbf{\textit{Equilibrium point}} ($E_x^{[6]},E_y^{[6]},E_z^{[6]}$): Value of $E_x^{[6]}$ is given by the positive real root of the equation given by,\\
$X^8 - 200X^7 + 9770.73X^6 + 23631.746X^5 - 56491.68254X^4 - 130234.92X^3 + 110210X^2 + 154819X - 121904.7619 = 0$\\
which means equilibrium point is $(1.951629, 0.470157934, 0.510619637)$. Here, $x^* = 1.951629$, which is only positive real root of equation above and $z^* = 0.510619637 < \min\{10,98.0778\}$, that means existence condition is satisfied for the equilibrium point in which the coexistence of predator, prey and scavenger are present. Putting these values, the jacobian matrix corresponding to this equilibrium point turns out to be,
\begin{center}
    $$J_{xyz} = 
\begin{bmatrix}
     -0.022 &   -0.98 & -0.98\\
   0.24 & 0 & 0.042\\
   0.26 &   0.027 & -0.02\\
\end{bmatrix}$$
\end{center}
For this matrix characteristic polynomial becomes, 
$$\lambda^3 + 0.041538\lambda^2 + 0.4892535\lambda + 0.02166$$ 
For this characteristic polynomial, $m_1 = 0.041538, m_2 = 0.4892535, m_3 = 0.02166$ which means $m_1, m_2, m_3 > 0$ and $m_1m_2 - m_3 = -0.001337 < 0$. That means the equilibrium point $(1.95, 0.47, 0.51)$ is unstable, but all three populations can coexist in the system.\\
As it is clear from Fig. \ref{fig:ex3} that all the populations coexist in the system as when the population of prey decreases, the food for scavengers and predators decreases. Hence, the population of predators and scavengers also decreases. But as the predator and scavengers can also rely on each other, their population decreases slower due to interrelation. As prey population decreases, scavenger population decreases, hence the predator population. As the predator population decreases the scavenger population starts decreasing much faster. As the population of predators and scavengers starts decreasing, prey population starts increasing.\\
\textbf{\textit{Remarks:}} Results for Example \ref{ex3} are as follows,
\begin{enumerate}
    \item[\textbf{1.}] Predator prey can coexist in the system in the absence of scavengers (Unstable).
    \vspace{-7.5pt}
    \item[\textbf{2.}] Scavenger prey can coexist in the system in the absence of predators (Unstable).
    \vspace{-7.5pt}
    \item[\textbf{3.}] Predator scavengers can't coexist in the system in the absence of prey.
    \vspace{-7.5pt}
    \item[\textbf{4.}] Predator prey scavengers can all coexist in the system (Unstable).
\end{enumerate}


\begin{example}\label{ex4}
    Consider predator prey scavenger system with parameter values such as $r = 1$, $k = 2$, $a = 1$, $a_0 = 0.25$, $b = 1$, $b_0 = 0.25$, $d = 1$, $e = 1$, $f = 1$, $g = 1$, $h = 0.25$, $i = 1$, $i_0 = 0.25$ and $j=1$, with initial conditions $x(0) = 4$, $y(0) = 3$ and $z(0) = 2$.
    \begin{equation*}
        \begin{split}
            \frac{dx}{dt} &= x\left(1 - \frac{x}{2}\right) - \frac{x^2 y}{1 + 0.25 x^2} - \frac{x^2 z}{1 + 0.25 x^2}\\
            \frac{dy}{dt} &= \frac{x^2 y}{1 + 0.25 x^2} + \frac{z^2 y}{1 + 0.25 z^2} - y\\
            \frac{dz}{dt} &= \frac{x^2 z}{1 + 0.25 x^2} + 0.25yz - \frac{yz^2}{1 + 0.25 z^2} - z
        \end{split}
    \end{equation*}
\end{example}

\begin{figure}[!h]
    \centering
    \begin{subfigure}[b]{7.5cm}
         \centering
         \includegraphics[width=\textwidth]{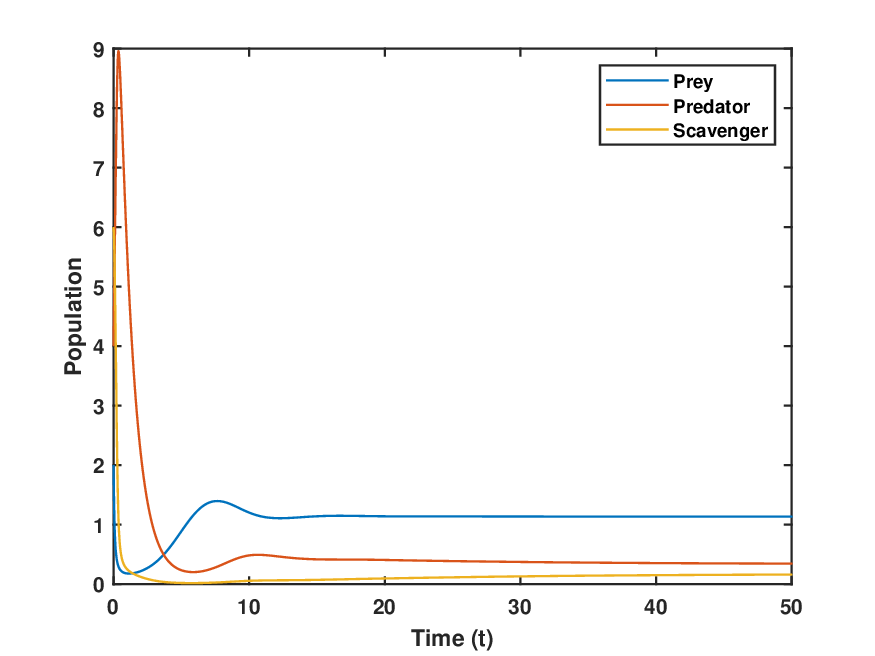}
         \caption{Dynamical behaviour.}
    \end{subfigure}
    \hfill
    \begin{subfigure}[b]{7.5cm}
         \centering
         \includegraphics[width=\textwidth]{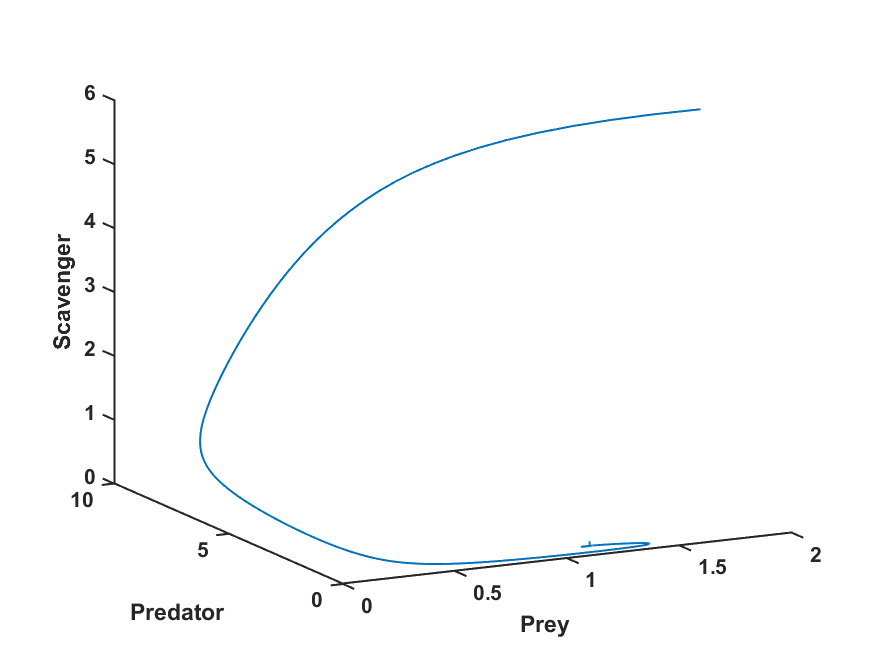}
         \caption{Phase portrait.}
    \end{subfigure}
        \caption{Predator prey scavenger model for Example \ref{ex4}.}
        \label{fig:ex4}
\end{figure}

For predator prey scavenger system, it is clear from Fig. \ref{fig:ex4} that all of the population gets stabilize at specific population values.
\\
\textbf{\textit{Equilibrium point}} ($0,0,0$): It is clear that $(0,0,0)$ is an unstable equilibrium point because if this point becomes stable, then the system collapse. As described in subsection \ref{PredPreyScav}, ($0,E_y,0$) and ($0,0,E_z$) cases are also included in this ($0,0,0$) case as both yields equilibrium point $(0,0,0)$ and both are not stable as well.
\\
\textbf{\textit{Equilibrium point}} ($E_x^{[2]},0,0$): For equilibrium point $(2,0,0)$, 
$$\frac{e}{d-ea_0} = \frac{4}{3} < 2^2 = k^2 \quad \text{ and } \quad \frac{j}{g-jb_0} = \frac{4}{3} < 2^2 = k^2$$ 
Hence, equilibrium point $(2,0,0)$ is not stable.
\\
\textbf{\textit{Equilibrium point}} ($0,E_y^{[3]},E_z^{[3]}$): As, $f-i_0 e = 0.75 > 0$ and
$$hf\left[ \frac{e}{f-i_0 e} \right]^{1/2} - ie = \frac{1-2\sqrt{3}}{2\sqrt{3}} = -0.7113248654 < 0$$ 
implies the existence conditions are not satisfied for the equilibrium point $\left(0,\frac{fj z_0}{hfz_0-ie}, z_0\right)$ as the equilibrium point turns out to have a negative number, which has no significance in nature. Hence, equilibrium point ($0,E_y^{[3]},E_z^{[3]}$) does not exist in the system.
\\
\textbf{\textit{Equilibrium point}} ($E_x^{[4]},E_y^{[4]},0$): As, $d-a_0 e = 0.75 > 0$ implies the existence condition for the equilibrium point, 
$$\left( \left[ \frac{e}{d-a_0 e} \right]^{1/2} , \frac{drx_0(k-x_0)}{aek}, 0\right) = \left(\frac{2}{\sqrt{3}},\frac{2(\sqrt{3}-1)}{3},0\right) = (1.1547,0.4880,0)$$
For this equilibrium point, 
$$\frac{2a_0 e(k-x_0)}{dk} = 0.211325 < 1$$ 
$$hrdx_0\left(1-\frac{x_0}{k}\right)[d+(b_0-a_0)e] + ae[ge - j[d+(b_0-a_0)e]] = 0.122 > 0$$ 
Hence, the stability condition are not met which means predator and prey may coexist in an unstable system.
\\
\textbf{\textit{Equilibrium point}} ($E_x^{[5]},0,E_z^{[5]}$): As, $g-b_0 j = 0.75 > 0$ implies the existence condition for the equilibrium point 
$$\left( \left[ \frac{j}{g-b_0 j} \right]^{1/2} , 0 , \frac{grx_0(k-x_0)}{bjk}\right) = \left(\frac{2}{\sqrt{3}},0,\frac{2(\sqrt{3}-1)}{3}\right) = (1.1547,0,0.4880)$$
For this equilibrium point, 
$$\frac{2b_0 j(k-x_0)}{gk} = 0.211325 < 1$$ 
$$t = \frac{e(g+j(a_0-b_0))-dj}{g+j(a_0-b_0)} = 0$$ 
and hence $z_0 ^2 = 0.238177 > \frac{t}{f-i_0 t} = 0$. Therefore, the stability condition doesn't satisfy that means scavenger and prey may coexist in an unstable system.
\\
\textbf{\textit{Equilibrium point}} ($E_x^{[6]},E_y^{[6]},E_z^{[6]}$): Value of $E_x^{[6]}$ is given by the positive real root of the equation given by,\\
$X^8 - 4X^7 + 3.926X^6 - 4.221845X^5 + 7.906X^4 + 6.77647X^3 - 5.163X^2 + 3.87227X - 17.19832 = 0$\\
which means equilibrium point is $(1.1331137,0.33726859,0.168040476)$. Here, $x^* = 1.1331137$, which is only positive real root of the equation above and $z^* = 0.168040476 < \min\{4,1.01062\}$, that means the existence condition is satisfied for the equilibrium point in which the coexistence of predator, prey and scavenger are present. Putting these values, the jacobian matrix corresponding to this equilibrium point turns out to be,
\begin{center}
    $$J_{xyz} = 
\begin{bmatrix}
     -0.79 &   -0.97 & -0.97\\
   0.44 & 0 & 0.11\\
   0.22 &   0.014 & -0.055\\
\end{bmatrix}$$
\end{center}
For this matrix characteristic polynomial becomes, 
$$\lambda^3 + 0.8448437\lambda^2 + 0.68008\lambda + 0.052045$$ 
For this characteristic polynomial, $m_1 = 0.845, m_2 = 0.68, m_3 = 0.052$ which means $m_1, m_2, m_3 > 0$ and $m_1m_2 - m_3 = 0.5225 > 0$. That means the equilibrium point $(1.133,0.337,0.168)$ is stable. Hence, all three populations can coexist in the system. As it is clear from Fig. \ref{fig:ex4} that all the populations coexist in the system. When the population of prey decreases, the food for scavengers and predators decreases, and hence the population of predators and scavengers also decreases. As the prey population decreases, the scavenger population also decreases and, hence the predator population. As the predator population decreases the scavenger population starts decreasing much faster. As the population of predators and scavengers starts decreasing, prey population starts increasing. Hence, the system gets stabilize, and all the population (prey, predator and scavengers) coexist in this system.\\ 
\textbf{\textit{Remarks:}} Results for Example \ref{ex4} are as follows,
\begin{enumerate}
    \item[\textbf{1.}] Predator prey can coexist in the system in the absence of scavengers (Unstable).
    \vspace{-7.5pt}
    \item[\textbf{2.}] Scavenger prey can coexist in the system in the absence of predators (Unstable).
    \vspace{-7.5pt}
    \item[\textbf{3.}] Predator scavengers can't coexist in the system in the absence of prey.
    \vspace{-7.5pt}
    \item[\textbf{4.}] Predator prey scavenger can all coexist in the system (Stable).
\end{enumerate}


\section{Parameter Estimation}\label{param_estimation}
Parameter estimation require proper dataset for the formulation and preparation of the model for the same. In this study population of various species from the American wildlife forest dataset \cite{nndata}. Two step model for parameter estimation is done in this study with first step consists of physics-informed deep neural network using adam backpropagation defined in section \ref{adam} and second step consist of fine-tune optimization of parameters using BFGS optimization algorihtm defined in section \ref{BFGS}.

\subsection{Problem Formulation}\label{problem_formulation}
\begin{figure}[!h]
    \centering
    \includegraphics[width=0.6\textwidth]{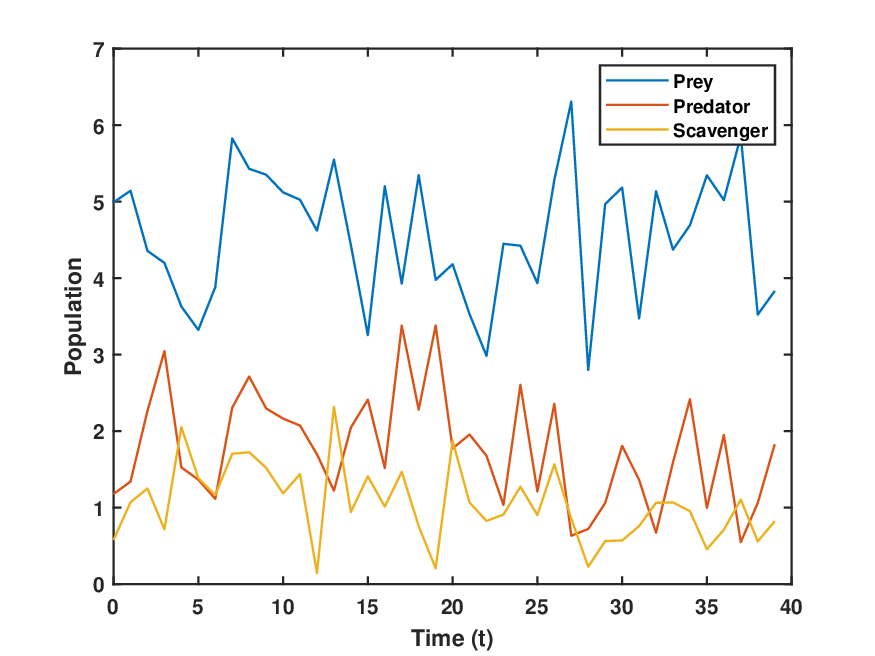}
    \caption{American forest dataset \cite{nndata}.}
    \label{fig:dataset}
\end{figure}

Fig. \ref{fig:dataset} represents the dataset used in this study. Dataset consists of population of various species across 40 years (1970-2009). Hence, initial assumptions for categorization of population species in the dataset and parameters of the proposed model \eqref{eq:propose} are as follows,
\begin{enumerate}
    \item[\textbf{1.}] Apex predators are assumed as predators.
    \vspace{-7.5pt}
    \item[\textbf{2.}] Scavengers are rotten dead body eater species.
    \vspace{-7.5pt}
    \item[\textbf{3.}] Rest of the species are taken as prey.
    \vspace{-7.5pt}
    \item[\textbf{4.}] All the parameters must have positive values.
\end{enumerate}

After applying these assumptions to the dataset, the dataset is then normalized using min-max scaling over independent parameter time (years) and numerical population of different species (predator, prey and scavenger).

Consider a dynamical system of $n$ ordinary differential equation with $m$ unknown parameters,
\begin{equation}\label{eq:dynamicalsystem}
\begin{split}
    \frac{dX(t)}{dt} &= f(t,X(t),p), \qquad t \in [t_0 , t_{end}] \subset \mathcal{R},\\
    X(t_0) &= X_0, \qquad X_0 \in \mathcal{R}^n,
\end{split}
\end{equation}
\par
where $X(t)= [x_1(t),x_2(t),\dots,x_n(t)]^{T} \in \mathcal{R}^n$ is population vector, $p = [p_1,p_2,\dots,p_m]\in \mathcal{R}^m$ is parameter vector, and $f$ is a vector function defined by,
\begin{center}
$f(t,X(t),p)=
    \begin{bmatrix}
      f_1 (t,x_1,x_2,\dots,x_n,p_1,p_2,\dots,p_m)\\
      f_2 (t,x_1,x_2,\dots,x_n,p_1,p_2,\dots,p_m)\\
      \dots\\
      f_n (t,x_1,x_2,\dots,x_n,p_1,p_2,\dots,p_m)
    \end{bmatrix}$
\end{center}
\par
Now in proposed model \eqref{eq:propose}, number of population variables is $n=3$ and number of parameters is $m=14$. Here, the problem becomes finding the unknown vector $p$ using other factors such as, physics-informed error (PIE), mean square error (MSE), and given data which will be explained in next subsections. 

\subsection{Physics-Informed Deep Neural Network Model}\label{PIDNN_model}
\begin{figure}[!h]
    \centering
    \includegraphics[width=\textwidth]{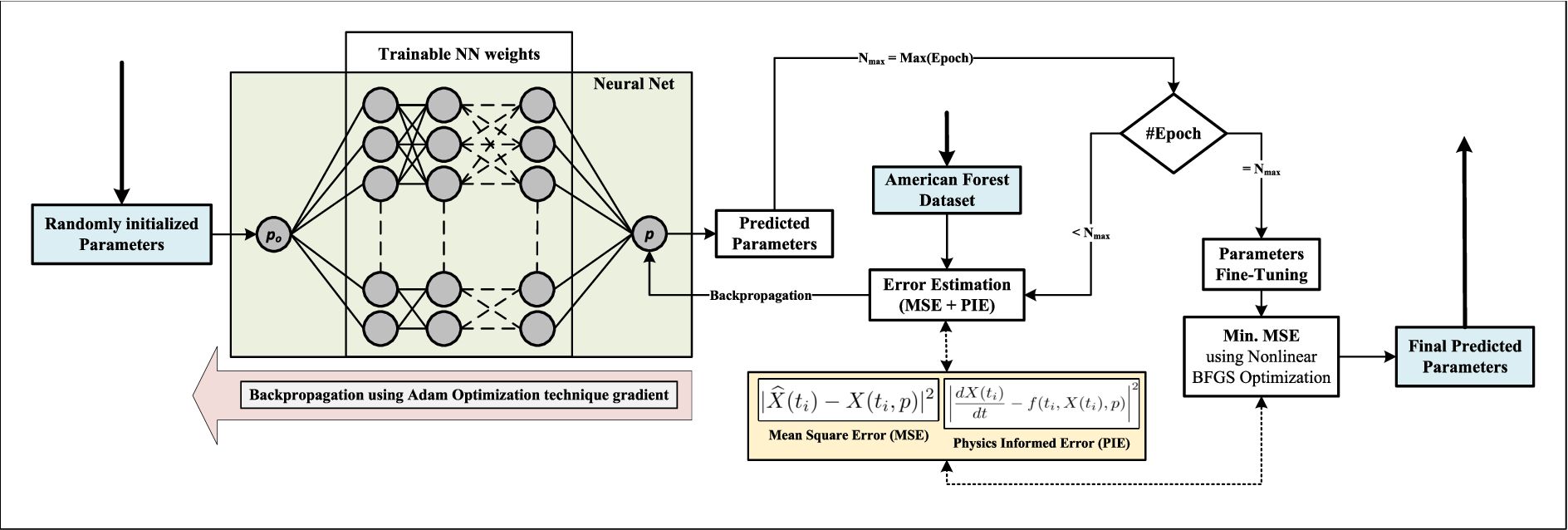}
    \caption{Flowchart for physics-informed deep neural network model.}
    \label{fig:flowchart}
\end{figure}

Fig. \ref{fig:flowchart} represents the physics-informed deep neural network model used for estimating the parameters of a dynamical system \eqref{eq:propose}. First, initialize the parameters of the dynamical system normally. Then random parameters are input into the neural network with the Adam backpropagation algorithm, and the objective function consists of MSE and PIE. Initialize the optimization process of neural network and update the weights and biases up to specific epoch limits using the dataset provided. The hence-found parameters are then input in the BFGS optimization algorithm for fine-tuning the parameter values using the objective function as MSE. The hence-found parameter values are the final predicted values for the parameters using the dataset.
\par
The physics-informed neural network (PINN) uses physical laws of differential equations into consideration of the loss function. Initialize the parameters $p_{init} = [p_1,p_2,\dots,p_{14}]\in \mathcal{R}^{14}$ using lognormal distribution that is,

\begin{equation*}
    \text{log}\left(p_i\right) \sim \mathcal{N}(\mu,\,\sigma^{2})\,\qquad \text{where, } i\in[1,2,\dots,14], \qquad \text{with } \mu = 0, \text{ and } \sigma = 1.
\end{equation*}

Considering a neural network with a backpropagation algorithm is adapted as Adam algorithm defined in subsection \ref{adam}, the required neural network can be presented as,
\begin{equation*}
    NN^{\Theta} : \mathcal{R}^{14} \rightarrow \mathcal{R}^{14} \qquad \text{where, trainable parameters, }\Theta = \{W,b\}
\end{equation*}

where, $\{W,b\}$ is a collection of trainable parameters, namely weights, and biases. The neural network defined above is trained to minimize the combined error (Mean Square Error and Physics-Informed Error) and predicts the parameters. Now, the objective function or, loss function (for minimization) for which PINN optimization is required can be written as,
\begin{equation}
    arg \min(MSE + PIE)
\end{equation}

where MSE is the mean square error term and PIE is the physics-informed error term. Let $p$ be the predicted parameter vector, $T_d = 40$ be the total data rows available, $\widehat{X}(t_i)$ be the true population values at time $t_i$ and $X(t_{i},p)$ be the predicted population value at time $t_i$ using predicted parameters $p$. Now, the overall error term can be represented as, 
\begin{equation}\label{eq:fullerror}
    E(p) = \frac{1}{|T_d|} \sum_{i\in T_d}  \left(\left|\widehat{X}(t_i) - X(t_{i},p)\right|^2 + \left|\frac{dX(t_i)}{dt}-f(t_i,X(t_i),p)\right|^2\right)
\end{equation}

where, $\left|\widehat{X}(t_i) - X(t_{i},p)\right|^2$ is the mean square error and $\left|\frac{dX(t_i)}{dt}-f(t_i,X(t_i),p)\right|^2$ is the physics-informed error. Combined error is normalized using the total number of data values $T_d$.
\par
The neural network structure used in this research is as follows, the input layer requires the input shape of vector to be $(14,)$, then $3$ hidden layers are implemented, each with $32$ neurons along with swish activation function. There are $14$ neurons in the output layer with Rectified Linear Unit (ReLU) activation function and the output shape of the predicted parameter vector is the same as input vector $(14,)$. The maximum epoch limit for the neural network is set to $100$ epochs.
\par
After the neural network, the predicted parameters are then fine-tuned using the BFGS optimization algorithm defined in subsection \ref{BFGS} with objective function as,
\begin{equation}
    arg \min(MSE)
\end{equation}

where MSE is defined the same as in equation \eqref{eq:fullerror}. That is, $\text{MSE} = \left|\widehat{X}(t_i) - X(t_{i},p)\right|^2$. The maximum iteration limit for the BFGS optimization algorithm is $200$ and threshold movement is also defined for BFGS optimization in physics physics-informed deep neural network model given by Fig. \ref{fig:flowchart}. Resulted parameters are final predicted parameters obtained from the dataset.

\subsection{Error Analysis and Results}\label{error_analysis}
\begin{figure}[!h]
    \centering
    \includegraphics[width=0.6\textwidth]{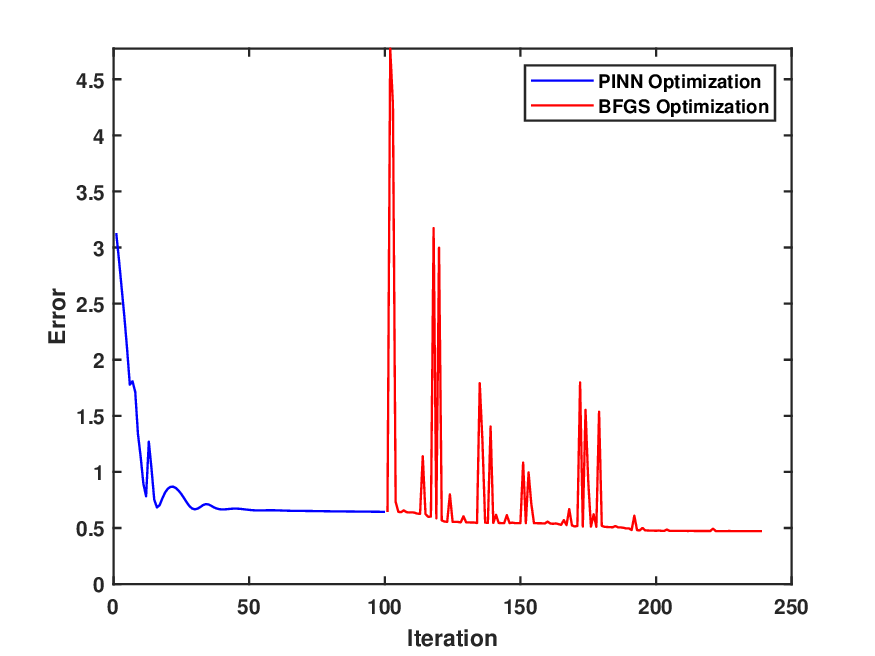}
    \caption{Error analysis: Total error propagation.}
    \label{fig:error}
\end{figure}

Fig. \ref{fig:error} represents the total error (MSE $+$ PIE) propagation over iterations. Physics-informed neural network with Adam backpropagation runs for $100$ epochs with initial total error $3.128$ which propagates through epochs and final error after neural net optimizes to $0.644$. After the neural net, the parameters are input into the BFGS optimization algorithm for fine-tuning. Now, the initial total error for BFGS $0.644$ optimizes to $0.472$. The final parameters are predicted parameters for model \eqref{eq:propose} using a natural dataset.

After physics-informed deep neural network, the final predicted parameter for model \eqref{eq:propose} is given by Table \ref{tab:real_values}. These parameter values define the underlying dynamics of the model \eqref{eq:propose}. Value $r$ defines that the population of prey increases with a decreasing rate in nature where total resources for prey is $k$, which is very large as compared to population values. Values $a_0$, $b_0$, and $i_0$ are handling time required by respective predators for their respective prey. It can be easily concluded that handling time is the least required by scavengers for handling prey populations. Values $d$, $f$ and $g$ are the rate of change of population of predators and scavengers and $h$ is the scavenge factor for the scavengers. From here it can be concluded that scavengers rely on the dead bodies of predators more than the prey population.
\begin{table}[!h]
    \centering
    \begin{tabular}{|c|c|c|} \hline 
         \textbf{S. No.}&  \textbf{Parameter}& \textbf{Numerical values}\\ \hline 
         01.& $r$ & $0.9701107742719246$\\ \hline
         02.& $k$ & $573.2545487545212544$\\ \hline
         03.& $a$ & $0.7668876233328743$\\ \hline
         04.& $a_0$ & $0.4686878655732233$\\ \hline
         05.& $b$ & $0.6893067418603573$\\ \hline
         06.& $b_0$ & $0.053266947840986595$\\ \hline
         07.& $d$ & $0.42441058569930494$\\ \hline
         08.& $f$ & $0.46630691773424437$\\ \hline
         09.& $e$ & $0.888598932493589$\\ \hline
         10.& $g$ & $0.08334616995047056$\\ \hline
         11.& $h$ & $0.16502232050920586$\\ \hline
         12.& $i$ & $1.05992612257741696$\\ \hline
         13.& $i_0$ & $0.105259076974745925$\\ \hline
         14.& $j$ & $0.5320956432008955$\\ \hline
    \end{tabular}
    \caption{Resulting values of parameters}
    \label{tab:real_values}
\end{table}

\begin{figure}[!h]
    \centering
    \begin{subfigure}[b]{7.5cm}
         \centering
         \includegraphics[width=\textwidth]{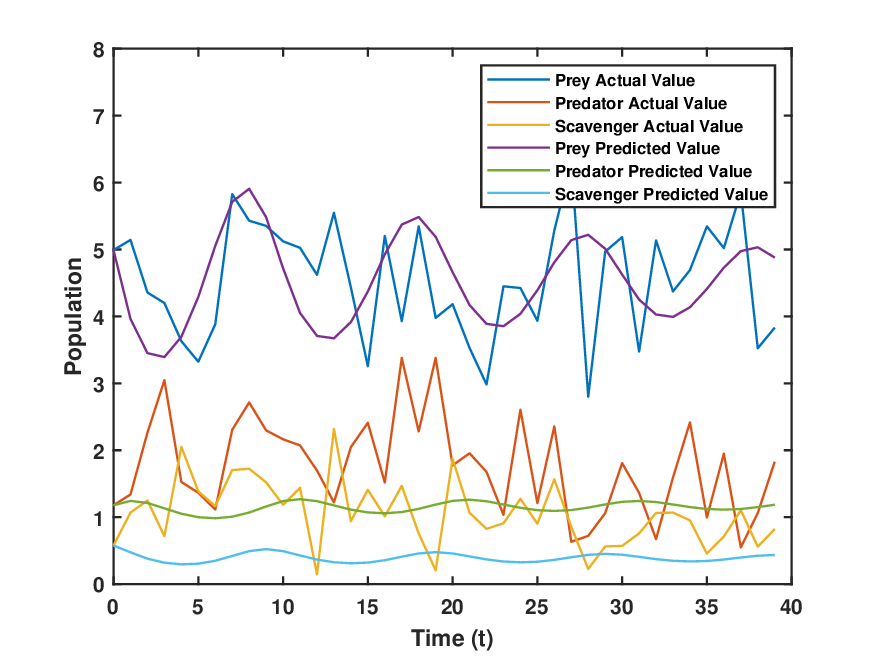}
         \caption{Comparison Graph.}
         \label{fig:parametersa}
    \end{subfigure}
    \hfill
    \begin{subfigure}[b]{7.5cm}
         \centering
         \includegraphics[width=\textwidth]{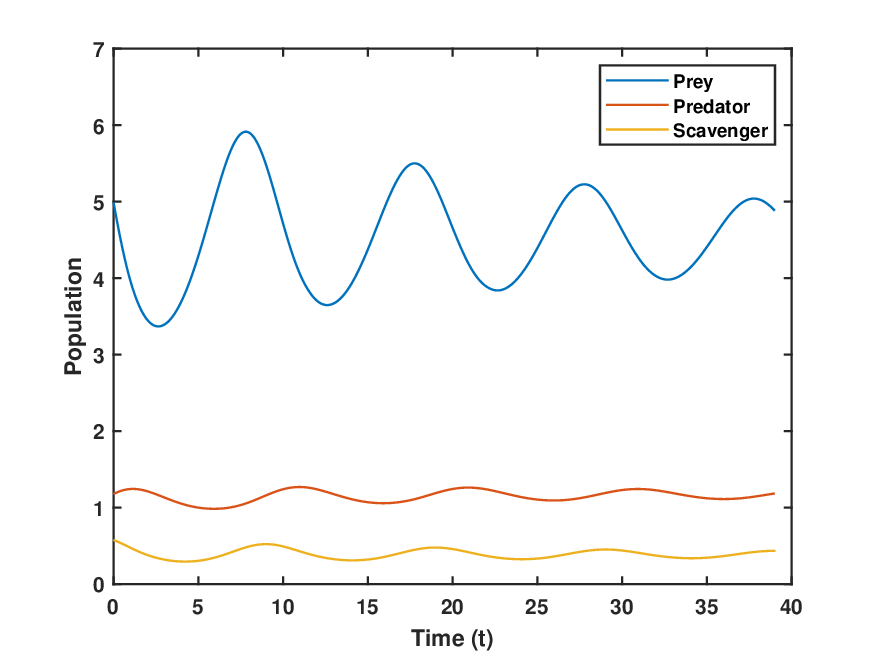}
         \caption{Plot for predicted parameters.}
         \label{fig:parametersb}
    \end{subfigure}
        \caption{Parameter accuracy comparison analysis.}
        \label{fig:parameters}
\end{figure}

Fig. \ref{fig:parametersa} represents the comparison graph of the normalized natural dataset with the model \eqref{eq:propose} using final predicted parameters. Fig. \ref{fig:parametersb} represents the plot for the predator-prey scavenger model using the parameters given in Table \ref{tab:real_values}. The comparison graph shown in Fig. \ref{fig:parametersa} shows that predicted parameters using a physics-informed deep neural network model explained in section \ref{PIDNN_model} captures most of the features of the natural dataset and mimics the data appropriately. Further, these predicted parameters can be more optimized and fine-tuned when the large dataset is available as the present dataset consists of only $40$ data point that is, population values.   

\section{Stability Analysis of Parameters}\label{parameter_analysis}
Using parameters obtained from Table \ref{tab:real_values} with initial values as $x(0)=4.991$, $y(0)=1.178$ and $z(0)=0.577$, the predator prey scavenger model \eqref{eq:propose} takes the form,
\begin{equation*}
    \begin{split}
        \frac{dx}{dt} &= 0.97011x\left(1-\frac{x}{573.25455}\right) - \frac{0.76689x^2 y}{1 + 0.46869 x^2} - \frac{0.68930x^2 z}{1 + 0.05327 x^2}\\
        \frac{dy}{dt} &= \frac{0.42441x^2 y}{1 + 0.46869 x^2} + \frac{0.46630z^2 y}{1 + 0.10526 z^2} - 0.88860y\\
        \frac{dz}{dt} &= \frac{0.08335x^2 z}{1 + 0.05327 x^2} + 0.16502yz - \frac{1.05992yz^2}{1 + 0.10526 z^2} - 0.53210z
    \end{split}
\end{equation*}

\begin{figure}[!h]
    \centering
    \begin{subfigure}[b]{7.5cm}
         \centering
         \includegraphics[width=\textwidth]{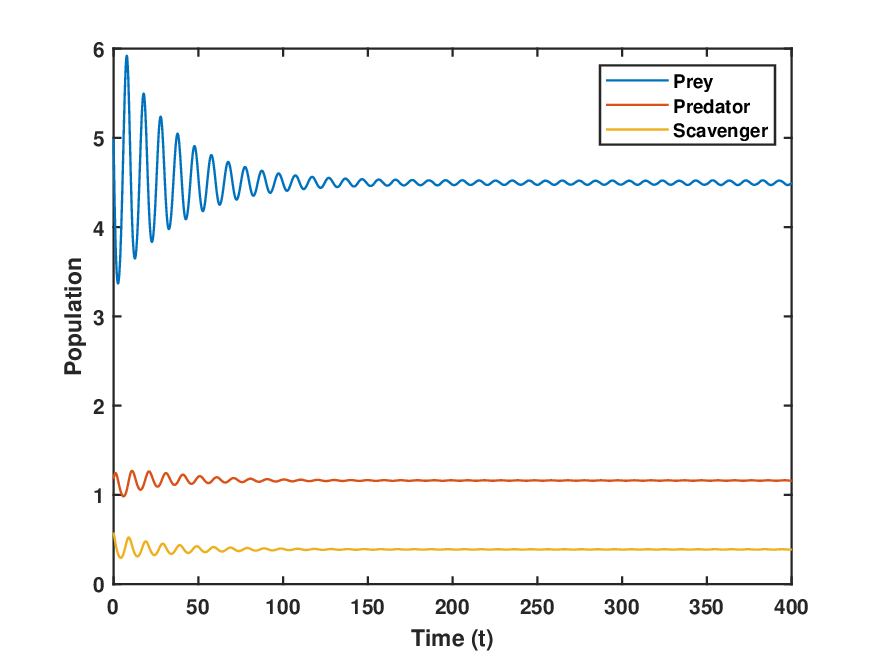}
         \caption{Dynamical behaviour.}
    \end{subfigure}
    \hfill
    \begin{subfigure}[b]{7.5cm}
         \centering
         \includegraphics[width=\textwidth]{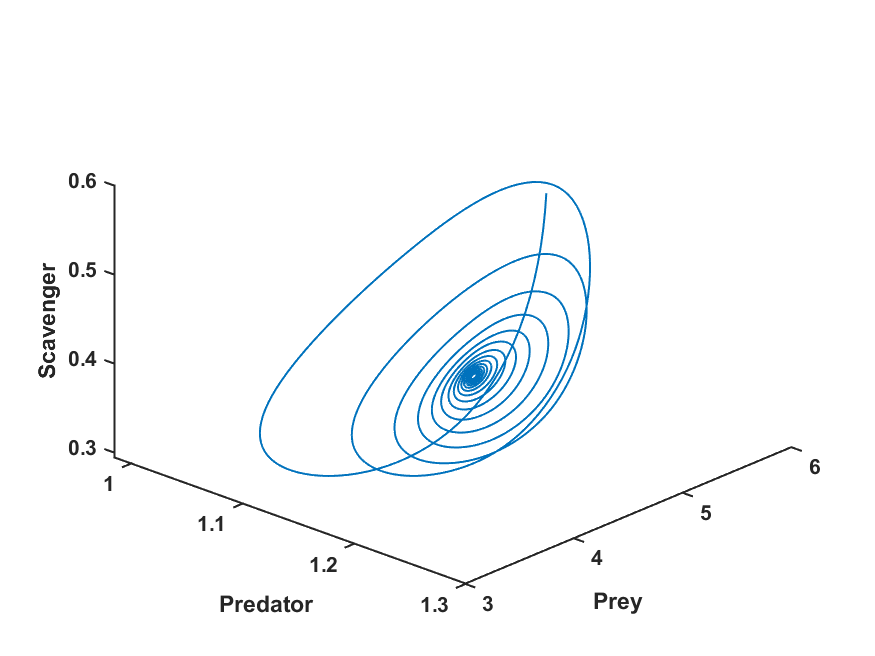}
         \caption{Phase portrait.}
    \end{subfigure}
        \caption{Predator prey scavenger model using parameters obtained in section \ref{param_estimation}.}
        \label{fig:real}
\end{figure}

For predator prey scavenger system using parameters obtained in Table \ref{tab:real_values}, it is clear from Fig. \ref{fig:real} that all of the population oscillates about specific population values.
\\
\textbf{\textit{Equilibrium point}} ($0,0,0$): It is clear that $(0,0,0)$ is an unstable equilibrium point because if this point becomes stable, then the system collapse. As described in subsection \ref{PredPreyScav}, ($0,E_y,0$) and ($0,0,E_z$) cases are also included in this ($0,0,0$) case as both yields equilibrium point $(0,0,0)$ and both are not stable as well.
\\
\textbf{\textit{Equilibrium point}} ($E_x^{[2]},0,0$): For equilibrium point $(573.25455,0,0)$, 
$$\frac{e}{d-ea_0} = 111.984 < (573.25455)^2 = k^2 \quad \text{ and } \quad \frac{j}{g-jb_0} = 9.673966 < (573.25455)^2 = k^2$$ 
Hence, equilibrium point $(573.25455,0,0)$ is not stable.
\\
\textbf{\textit{Equilibrium point}} ($0,E_y^{[3]},E_z^{[3]}$): As, $f-i_0 e = 0.3727738 > 0$ and
$$hf\left[ \frac{e}{f-i_0 e} \right]^{1/2} - ie = -0.8230414 < 0$$ 
implies the existence conditions are not satisfied for the equilibrium point $\left(0,\frac{fj z_0}{hfz_0-ie}, z_0\right)$ as the equilibrium point turns out to have a negative number, which has no significance in nature. Hence, equilibrium point ($0,E_y^{[3]},E_z^{[3]}$) does not exist in the system.
\\
\textbf{\textit{Equilibrium point}} ($E_x^{[4]},E_y^{[4]},0$): As, $d-a_0 e = 0.007935 > 0$ implies the existence condition for the equilibrium point, 
$$\left( \left[ \frac{e}{d-a_0 e} \right]^{1/2} , \frac{drx_0(k-x_0)}{aek}, 0\right) = (10.582252,6.2756,0)$$
For this equilibrium point, 
$$\frac{2a_0 e(k-x_0)}{dk} = 1.9263771 > 1$$
Hence, the stability condition are not met which means predator and prey may coexist in an unstable system.
\\
\textbf{\textit{Equilibrium point}} ($E_x^{[5]},0,E_z^{[5]}$): As, $g-b_0 j = 0.055003 > 0$ implies the existence condition for the equilibrium point 
$$\left( \left[ \frac{j}{g-b_0 j} \right]^{1/2} , 0 , \frac{grx_0(k-x_0)}{bjk}\right) = (3.1103,0,0.6819357)$$
For this equilibrium point, 
$$\frac{2b_0 j(k-x_0)}{gk} = 0.67643973 < 1$$ 
$$t = \frac{e(g+j(a_0-b_0))-dj}{g+j(a_0-b_0)} = 0.1467$$ 
and hence $z_0 ^2 = 0.465 > \frac{t}{f-i_0 t} = 0.32537$. Therefore, the stability condition doesn't satisfy that means scavenger and prey may coexist in an unstable system.
\\
\textbf{\textit{Equilibrium point}} ($E_x^{[6]},E_y^{[6]},E_z^{[6]}$): Value of $E_x^{[6]}$ is given by the positive real root of the equation \eqref{eq:appendix} using parameter values given in Table \ref{tab:real_values}, which means equilibrium point is $(4.4984538,1.161178,0.38895175)$. Here, $x^* = 4.4984538$, which is only positive real root of the equation above and $z^* = 0.38895175 < \min\{12.239589,369.74242836\}$, that means the existence condition is satisfied for the equilibrium point in which the coexistence of predator, prey and scavenger are present. Putting these values, the jacobian matrix corresponding to this equilibrium point turns out to be,
\begin{center}
    $$J_{xyz} = 
\begin{bmatrix}
     0.34 &   -1.5 & -6.7\\
   0.04 & 0 & 0.41\\
   0.068 &   -0.094 & -0.46\\
\end{bmatrix}$$
\end{center}
For this matrix characteristic polynomial becomes, 
$$\lambda^3 + 0.1178643784\lambda^2 + 0.39684137\lambda + 0.02975781$$ 
For this characteristic polynomial, $m_1 = 0.1178643784, m_2 = 0.39684137, m_3 = 0.02975781$ which means $m_1, m_2, m_3 > 0$ and $m_1m_2 - m_3 = 0.01701565148 > 0$. That means the equilibrium point $(4.4984538,1.161178,0.38895175)$ is stable. Hence, all three populations can coexist in the system. As it is clear from Fig. \ref{fig:real} all the populations coexist in the system. When the population of prey decreases, the food for scavengers and predators decreases, and hence the population of predators and scavengers also decreases. As the prey population decreases, the scavenger population also decreases and, hence the predator population. As the predator population decreases the scavenger population starts decreasing much faster. As the population of predators and scavengers starts decreasing, the prey population starts increasing. Hence, the system gets stabilized, and all the population (prey, predator, and scavengers) coexist in this system.\\
\textbf{\textit{Remarks:}} Results for hence found real parameters are as follows,
\begin{enumerate}
    \item[\textbf{1.}] Predator-prey can coexist in the system in the absence of scavengers (Unstable).
    \vspace{-7.5pt}
    \item[\textbf{2.}] Scavenger prey can coexist in the system in the absence of predators (Unstable).
    \vspace{-7.5pt}
    \item[\textbf{3.}] Predator scavengers can't coexist in the system in the absence of prey.
    \vspace{-7.5pt}
    \item[\textbf{4.}] Predator prey scavenger can all coexist in the system (Stable).
\end{enumerate}


\section{Conclusion and Discussion}\label{conclusion}
In conclusion, this research aims to address the dynamic behavior of predator prey scavenger interactions in natural ecosystems through the development and analysis of a mathematical model. By extending classical Lotka-Volterra dynamics to incorporate scavenger behavior, this study aims to enhance the understanding of complex interactions between different species present in the ecosystem. The proposed model, based on principles of population dynamics and Holling's functional response of type III, provides a framework for simulating the interactions between predators, prey, and scavenger populations under different scenarios. Through numerical simulations, the model demonstrates robustness and versatility, offering insights into non-linear dynamics present in ecological systems. Real-world validation is achieved through parameter estimation using a natural population dataset provided by the Forest Service, U.S. Department. Leveraging computational techniques, physics-informed neural networks using Adam backpropagation and BFGS optimization algorithm, the model's parameters are predicted which approximately reflect observed ecological patterns. Stability analysis, conducted via Jacobian matrix analysis, ensures the reliability of the model under different conditions. By examining the stability of steady states, this analysis provides crucial insights into the long-term behavior of different populations namely, predator, prey, and scavenger present in the ecosystem.
\par
The proposed predator prey scavenger model using Holling’s functional response of type III of order $n=2$ is used for respective predator-prey interactions. Two-step parameter estimation from the natural dataset is also done in this paper using physics-informed deep neural network and BFGS optimization for further fine-tuning the estimated parameters. Numerical simulations and stability analysis of the proposed predator-prey scavenger model are also done for different scenarios. Looking ahead, this research sets the stage for further exploration and refinement of the model, by minimizing parameter-induced errors, conducting bifurcation analysis, and sensitivity analysis of the parameters to enhance predictive accuracy. Incorporating scavenger population and scavenge factor into mathematical modeling, this research aims to bridge the gap between theoretical modeling and actual population. Hence, this study also contributes to the advancement of ecological science and the sustainable management of natural resources.

\section*{Ethics declarations}
\subsection*{Conflict of interest}
The authors declare no conflict of interest.

\subsection*{Ethical approval}
This article does not contain any studies with human participants and/or animals performed by the author.

\subsection*{Data availability}
Wildlife population in U.S. dataset, is available to download from Forest Service, U.S. Dept. of Agriculture (USDA) repository at \href{https://doi.org/10.2737/RMRS-GTR-296}{https://doi.org/10.2737/RMRS-GTR-296} \cite{nndata}

\subsection*{Funding}
This research received no specific grant from any funding agency in the public, commercial, or not-for-profit sectors.


\bibliographystyle{elsarticle-num}
\bibliography{main}

\begin{thebibliography}{10}
\expandafter\ifx\csname url\endcsname\relax
  \def\url#1{\texttt{#1}}\fi
\expandafter\ifx\csname urlprefix\endcsname\relax\def\urlprefix{URL }\fi
\expandafter\ifx\csname href\endcsname\relax
  \def\href#1#2{#2} \def\path#1{#1}\fi

\bibitem{applications}
M.~L. Abell, J.~Braselton, Modern differential equations: Theory, applications, technology (2001).

\bibitem{Preyrefugee}
S.~Sarwardi, P.~K. Mandal, S.~Ray, Dynamical behaviour of a two-predator model with prey refuge, Journal of biological physics 39~(4) (2013) 701--722.

\bibitem{Trophicchains}
R.~Dilao, T.~Domingos, A general approach to the modelling of trophic chains, Ecological Modelling 132~(3) (2000) 191--202.

\bibitem{diseaseprey}
S.~Kant, V.~Kumar, Stability analysis of predator--prey system with migrating prey and disease infection in both species, Applied Mathematical Modelling 42 (2017) 509--539.

\bibitem{popdynamics1}
H.~Shim, P.~Fishwick, Visualization and interaction design for ecosystem modeling (2008).

\bibitem{popdynamics2}
T.~D. Schowalter, Insect ecology: an ecosystem approach, Academic press, 2022.

\bibitem{popdynamics3}
A.~Hastings, L.~Gross, Encyclopedia of theoretical ecology, no.~4, Univ of California Press, 2012.

\bibitem{LotkaVolterra}
M.-C. Anisiu, Lotka, volterra and their model, Did{\'a}ctica mathematica 32~(01) (2014).

\bibitem{LotVol}
V.~Had{\v{z}}iabdi{\'c}, M.~Mehulji{\'c}, J.~Bekte{\v{s}}evi{\'c}, Lotka-volterra model with two predators and their prey, Tem Journal 6~(1) (2017) 2217--8309.

\bibitem{figures}
S.-X. Wu, X.-Y. Meng, et~al., Dynamics of a delayed predator-prey system with fear effect, herd behavior and disease in the susceptible prey, AIMS Mathematics 6~(4) (2021) 3654--3685.

\bibitem{delayholling}
C.~Arora, V.~Kumar, A delayed version of prey predator system with modified holling-tanner response, in: AIP Conference Proceedings, Vol. 2336, AIP Publishing, 2021.

\bibitem{ann}
I.~A. Basheer, M.~Hajmeer, Artificial neural networks: fundamentals, computing, design, and application, Journal of microbiological methods 43~(1) (2000) 3--31.

\bibitem{odeNNbook}
N.~Yadav, A.~Yadav, M.~Kumar, et~al., An introduction to neural network methods for differential equations, Vol.~1, Springer, 2015.

\bibitem{neuron}
W.~S. McCulloch, W.~Pitts, A logical calculus of the ideas immanent in nervous activity, The bulletin of mathematical biophysics 5 (1943) 115--133.

\bibitem{PINNpde}
M.~Raissi, P.~Perdikaris, G.~E. Karniadakis, Physics-informed neural networks: A deep learning framework for solving forward and inverse problems involving nonlinear partial differential equations, Journal of Computational physics 378 (2019) 686--707.

\bibitem{nndata}
C.~H. Flather, M.~S. Knowles, M.~F. Jones, C.~Schilli, Wildlife population and harvest trends in the united states: a technical document supporting the forest service 2010 rpa assessment, General Technical Report RMRS-GTR-296, U.S. Department of Agriculture, Forest Service, Rocky Mountain Research Station, Fort Collins, CO (2013).

\bibitem{hollingintro}
J.~Dawes, M.~Souza, A derivation of holling's type i, ii and iii functional responses in predator--prey systems, Journal of theoretical biology 327 (2013) 11--22.

\bibitem{Hol1}
Y.~Zhang, X.~Meng, Dynamic study of a stochastic holling iii predator-prey system with a prey refuge, IFAC-PapersOnLine 55~(3) (2022) 73--78.

\bibitem{Hol2}
Y.~Huang, F.~Chen, L.~Zhong, Stability analysis of a prey--predator model with holling type iii response function incorporating a prey refuge, Applied Mathematics and Computation 182~(1) (2006) 672--683.

\bibitem{Hol3}
R.~Banerjee, P.~Das, D.~Mukherjee, Stability and permanence of a discrete-time two-prey one-predator system with holling type-iii functional response, Chaos, Solitons \& Fractals 117 (2018) 240--248.

\bibitem{adam}
D.~P. Kingma, J.~Ba, Adam: A method for stochastic optimization, arXiv preprint arXiv:1412.6980 (2014).

\bibitem{bfgs-b}
C.~G. Broyden, The convergence of a class of double-rank minimization algorithms 1. general considerations, IMA Journal of Applied Mathematics 6~(1) (1970) 76--90.

\bibitem{bfgs-f}
R.~Fletcher, A new approach to variable metric algorithms, The computer journal 13~(3) (1970) 317--322.

\bibitem{bfgs-g}
D.~Goldfarb, A family of variable-metric methods derived by variational means, Mathematics of computation 24~(109) (1970) 23--26.

\bibitem{bfgs-s}
D.~F. Shanno, Conditioning of quasi-newton methods for function minimization, Mathematics of computation 24~(111) (1970) 647--656.

\bibitem{bfgs-book}
J.~Nocedal, S.~J. Wright, Numerical optimization, Springer, 1999.

\bibitem{routh-hur}
E.~X. DeJesus, C.~Kaufman, Routh-hurwitz criterion in the examination of eigenvalues of a system of nonlinear ordinary differential equations, Physical Review A 35~(12) (1987) 5288.

\end{thebibliography}

\newpage
\section*{Appendix: Coefficients of Equation \eqref{eq:appendix}}\label{appendix}
\footnotesize	{
\begin{enumerate}
    \item[$\mathbf{p_0} = $] $e^{3} i^{2} {i_0}^{2} k^{2} r^{2}-2 e^{2} f i^{2} {i_0} k^{2} r^{2}+e f^{2} h^{2} {i_0} k^{2} r^{2}+e f^{2} i^{2} k^{2} r^{2}-f^{3} h^{2} k^{2} r^{2}$

    \item[$\mathbf{p_1} = $] $-2 a e f^{2} h {i_0} j k^{2} r-2 e^{3} i^{2} {i_0}^{2} k r^{2}+2 a f^{3} h j k^{2} r+4 e^{2} f i^{2} {i_0} k r^{2} -2 e f^{2} h^{2} {i_0} k r^{2} -2 e f^{2} i^{2} k r^{2}+2 f^{3} h^{2} k r^{2}$

    \item[$\mathbf{p_2} = $] $3 {a_0} e^{3} i^{2} {i_0}^{2} k^{2} r^{2}+2 {b_0} e^{3} i^{2} {i_0}^{2} k^{2} r^{2}-6 {a_0} e^{2} f i^{2} {i_0} k^{2} r^{2}+3 {a_0} e f^{2} h^{2} {i_0} k^{2} r^{2} -4 {b_0} e^{2} f i^{2} {i_0} k^{2} r^{2}+2 {b_0} e f^{2} h^{2} {i_0} k^{2} r^{2}-3 d e^{2} i^{2} {i_0}^{2} k^{2} r^{2}+a^{2} e f^{2} {i_0} j^{2} k^{2} -2 a b e^{2} f i {i_0} j k^{2}+3 {a_0} e f^{2} i^{2} k^{2} r^{2}-3 {a_0} f^{3} h^{2} k^{2} r^{2}+b^{2} e^{3} i^{2} {i_0} k^{2}+2 {b_0} e f^{2} i^{2} k^{2} r^{2} -2 {b_0} f^{3} h^{2} k^{2} r^{2}+4 d e f i^{2} {i_0} k^{2} r^{2}-d f^{2} h^{2} {i_0} k^{2} r^{2}-a^{2} f^{3} j^{2} k^{2}+2 a b e f^{2} i j k^{2} +2 a e f^{2} h {i_0} j k r-b^{2} e^{2} f i^{2} k^{2}+b^{2} e f^{2} h^{2} k^{2}-d f^{2} i^{2} k^{2} r^{2}+e^{3} i^{2} {i_0}^{2} r^{2}-2 a f^{3} h j k r -2 e^{2} f i^{2} {i_0} r^{2}+e f^{2} h^{2} {i_0} r^{2}+e f^{2} i^{2} r^{2}-f^{3} h^{2} r^{2}$

    \item[$\mathbf{p_3} = $] $-4 a {a_0} e f^{2} h {i_0} j k^{2} r-4 a {b_0} e f^{2} h {i_0} j k^{2} r-6 {a_0} e^{3} i^{2} {i_0}^{2} k r^{2}-4 {b_0} e^{3} i^{2} {i_0}^{2} k r^{2} +4 a {a_0} f^{3} h j k^{2} r+4 a {b_0} f^{3} h j k^{2} r+2 a d f^{2} h {i_0} j k^{2} r+2 a e f^{2} g h {i_0} k^{2} r +12 {a_0} e^{2} f i^{2} {i_0} k r^{2}-6 {a_0} e f^{2} h^{2} {i_0} k r^{2}+8 {b_0} e^{2} f i^{2} {i_0} k r^{2}-4 {b_0} e f^{2} h^{2} {i_0} k r^{2} +6 d e^{2} i^{2} {i_0}^{2} k r^{2}-2 a f^{3} g h k^{2} r-6 {a_0} e f^{2} i^{2} k r^{2}+6 {a_0} f^{3} h^{2} k r^{2}-4 {b_0} e f^{2} i^{2} k r^{2} +4 {b_0} f^{3} h^{2} k r^{2}-8 d e f i^{2} {i_0} k r^{2}+2 d f^{2} h^{2} {i_0} k r^{2}+2 d f^{2} i^{2} k r^{2}$

    \item[$\mathbf{p_4} = $] $3 {a_0}^{2} e^{3} i^{2} {i_0}^{2} k^{2} r^{2}+6 {a_0} {b_0} e^{3} i^{2} {i_0}^{2} k^{2} r^{2}+{b_0}^{2} e^{3} i^{2} {i_0}^{2} k^{2} r^{2}-6 {a_0}^{2} e^{2} f i^{2} {i_0} k^{2} r^{2} +3 {a_0}^{2} e f^{2} h^{2} {i_0} k^{2} r^{2}-12 {a_0} {b_0} e^{2} f i^{2} {i_0} k^{2} r^{2}+6 {a_0} {b_0} e f^{2} h^{2} {i_0} k^{2} r^{2} -6 {a_0} d e^{2} i^{2} {i_0}^{2} k^{2} r^{2}-2 {b_0}^{2} e^{2} f i^{2} {i_0} k^{2} r^{2}+{b_0}^{2} e f^{2} h^{2} {i_0} k^{2} r^{2}-6 {b_0} d e^{2} i^{2} {i_0}^{2} k^{2} r^{2} +a^{2} {a_0} e f^{2} {i_0} j^{2} k^{2}+2 a^{2} {b_0} e f^{2} {i_0} j^{2} k^{2}-4 a {a_0} b e^{2} f i {i_0} j k^{2}-2 a b {b_0} e^{2} f i {i_0} j k^{2} +3 {a_0}^{2} e f^{2} i^{2} k^{2} r^{2}-3 {a_0}^{2} f^{3} h^{2} k^{2} r^{2}+3 {a_0} b^{2} e^{3} i^{2} {i_0} k^{2}+6 {a_0} {b_0} e f^{2} i^{2} k^{2} r^{2} -6 {a_0} {b_0} f^{3} h^{2} k^{2} r^{2}+8 {a_0} d e f i^{2} {i_0} k^{2} r^{2}-2 {a_0} d f^{2} h^{2} {i_0} k^{2} r^{2}+{b_0}^{2} e f^{2} i^{2} k^{2} r^{2} -{b_0}^{2} f^{3} h^{2} k^{2} r^{2}+8 {b_0} d e f i^{2} {i_0} k^{2} r^{2}-2 {b_0} d f^{2} h^{2} {i_0} k^{2} r^{2}+3 d^{2} e i^{2} {i_0}^{2} k^{2} r^{2} -a^{2} {a_0} f^{3} j^{2} k^{2}-2 a^{2} {b_0} f^{3} j^{2} k^{2}-a^{2} d f^{2} {i_0} j^{2} k^{2}-2 a^{2} e f^{2} g {i_0} j k^{2} +4 a {a_0} b e f^{2} i j k^{2}+4 a {a_0} e f^{2} h {i_0} j k r+2 a b {b_0} e f^{2} i j k^{2}+4 a b d e f i {i_0} j k^{2} +2 a b e^{2} f g i {i_0} k^{2}+4 a {b_0} e f^{2} h {i_0} j k r-3 {a_0} b^{2} e^{2} f i^{2} k^{2}+3 {a_0} b^{2} e f^{2} h^{2} k^{2} -2 {a_0} d f^{2} i^{2} k^{2} r^{2}+3 {a_0} e^{3} i^{2} {i_0}^{2} r^{2}-3 b^{2} d e^{2} i^{2} {i_0} k^{2}-2 {b_0} d f^{2} i^{2} k^{2} r^{2} +2 {b_0} e^{3} i^{2} {i_0}^{2} r^{2}-2 d^{2} f i^{2} {i_0} k^{2} r^{2}+2 a^{2} f^{3} g j k^{2}-4 a {a_0} f^{3} h j k r-2 a b d f^{2} i j k^{2} -2 a b e f^{2} g i k^{2}-4 a {b_0} f^{3} h j k r-2 a d f^{2} h {i_0} j k r-2 a e f^{2} g h {i_0} k r-6 {a_0} e^{2} f i^{2} {i_0} r^{2} +3 {a_0} e f^{2} h^{2} {i_0} r^{2}+2 b^{2} d e f i^{2} k^{2}-b^{2} d f^{2} h^{2} k^{2}-4 {b_0} e^{2} f i^{2} {i_0} r^{2}+2 {b_0} e f^{2} h^{2} {i_0} r^{2} -3 d e^{2} i^{2} {i_0}^{2} r^{2}+2 a f^{3} g h k r+3 {a_0} e f^{2} i^{2} r^{2}-3 {a_0} f^{3} h^{2} r^{2}+2 {b_0} e f^{2} i^{2} r^{2} -2 {b_0} f^{3} h^{2} r^{2}+4 d e f i^{2} {i_0} r^{2}-d f^{2} h^{2} {i_0} r^{2}-d f^{2} i^{2} r^{2}$

    \item[$\mathbf{p_5} = $] $-2 a {a_0}^{2} e f^{2} h {i_0} j k^{2} r-8 a {a_0} {b_0} e f^{2} h {i_0} j k^{2} r-2 a {b_0}^{2} e f^{2} h {i_0} j k^{2} r-6 {a_0}^{2} e^{3} i^{2} {i_0}^{2} k r^{2} -12 {a_0} {b_0} e^{3} i^{2} {i_0}^{2} k r^{2}-2 {b_0}^{2} e^{3} i^{2} {i_0}^{2} k r^{2}+2 a {a_0}^{2} f^{3} h j k^{2} r+8 a {a_0} {b_0} f^{3} h j k^{2} r +2 a {a_0} d f^{2} h {i_0} j k^{2} r+4 a {a_0} e f^{2} g h {i_0} k^{2} r+2 a {b_0}^{2} f^{3} h j k^{2} r+4 a {b_0} d f^{2} h {i_0} j k^{2} r +2 a {b_0} e f^{2} g h {i_0} k^{2} r+12 {a_0}^{2} e^{2} f i^{2} {i_0} k r^{2}-6 {a_0}^{2} e f^{2} h^{2} {i_0} k r^{2}+24 {a_0} {b_0} e^{2} f i^{2} {i_0} k r^{2} -12 {a_0} {b_0} e f^{2} h^{2} {i_0} k r^{2}+12 {a_0} d e^{2} i^{2} {i_0}^{2} k r^{2}+4 {b_0}^{2} e^{2} f i^{2} {i_0} k r^{2}-2 {b_0}^{2} e f^{2} h^{2} {i_0} k r^{2} +12 {b_0} d e^{2} i^{2} {i_0}^{2} k r^{2}-4 a {a_0} f^{3} g h k^{2} r-2 a {b_0} f^{3} g h k^{2} r-2 a d f^{2} g h {i_0} k^{2} r -6 {a_0}^{2} e f^{2} i^{2} k r^{2}+6 {a_0}^{2} f^{3} h^{2} k r^{2}-12 {a_0} {b_0} e f^{2} i^{2} k r^{2}+12 {a_0} {b_0} f^{3} h^{2} k r^{2} -16 {a_0} d e f i^{2} {i_0} k r^{2}+4 {a_0} d f^{2} h^{2} {i_0} k r^{2}-2 {b_0}^{2} e f^{2} i^{2} k r^{2}+2 {b_0}^{2} f^{3} h^{2} k r^{2} -16 {b_0} d e f i^{2} {i_0} k r^{2}+4 {b_0} d f^{2} h^{2} {i_0} k r^{2}-6 d^{2} e i^{2} {i_0}^{2} k r^{2}+4 {a_0} d f^{2} i^{2} k r^{2} +4 {b_0} d f^{2} i^{2} k r^{2}+4 d^{2} f i^{2} {i_0} k r^{2}$

    \item[$\mathbf{p_6} = $] ${a_0}^{3} e^{3} i^{2} {i_0}^{2} k^{2} r^{2} +6 {a_0}^{2} {b_0} e^{3} i^{2} {i_0}^{2} k^{2} r^{2}+3 {a_0} {b_0}^{2} e^{3} i^{2} {i_0}^{2} k^{2} r^{2}-2 {a_0}^{3} e^{2} f i^{2} {i_0} k^{2} r^{2} +{a_0}^{3} e f^{2} h^{2} {i_0} k^{2} r^{2}-12 {a_0}^{2} {b_0} e^{2} f i^{2} {i_0} k^{2} r^{2}+6 {a_0}^{2} {b_0} e f^{2} h^{2} {i_0} k^{2} r^{2} -3 {a_0}^{2} d e^{2} i^{2} {i_0}^{2} k^{2} r^{2}-6 {a_0} {b_0}^{2} e^{2} f i^{2} {i_0} k^{2} r^{2}+3 {a_0} {b_0}^{2} e f^{2} h^{2} {i_0} k^{2} r^{2} -12 {a_0} {b_0} d e^{2} i^{2} {i_0}^{2} k^{2} r^{2}-3 {b_0}^{2} d e^{2} i^{2} {i_0}^{2} k^{2} r^{2}+2 a^{2} {a_0} {b_0} e f^{2} {i_0} j^{2} k^{2} +a^{2} {b_0}^{2} e f^{2} {i_0} j^{2} k^{2}-2 a {a_0}^{2} b e^{2} f i {i_0} j k^{2}-4 a {a_0} b {b_0} e^{2} f i {i_0} j k^{2}+{a_0}^{3} e f^{2} i^{2} k^{2} r^{2} -{a_0}^{3} f^{3} h^{2} k^{2} r^{2}+3 {a_0}^{2} b^{2} e^{3} i^{2} {i_0} k^{2}+6 {a_0}^{2} {b_0} e f^{2} i^{2} k^{2} r^{2}-6 {a_0}^{2} {b_0} f^{3} h^{2} k^{2} r^{2} +4 {a_0}^{2} d e f i^{2} {i_0} k^{2} r^{2}-{a_0}^{2} d f^{2} h^{2} {i_0} k^{2} r^{2}+3 {a_0} {b_0}^{2} e f^{2} i^{2} k^{2} r^{2}-3 {a_0} {b_0}^{2} f^{3} h^{2} k^{2} r^{2} +16 {a_0} {b_0} d e f i^{2} {i_0} k^{2} r^{2}-4 {a_0} {b_0} d f^{2} h^{2} {i_0} k^{2} r^{2}+3 {a_0} d^{2} e i^{2} {i_0}^{2} k^{2} r^{2}+4 {b_0}^{2} d e f i^{2} {i_0} k^{2} r^{2} -{b_0}^{2} d f^{2} h^{2} {i_0} k^{2} r^{2}+6 {b_0} d^{2} e i^{2} {i_0}^{2} k^{2} r^{2}-2 a^{2} {a_0} {b_0} f^{3} j^{2} k^{2}-2 a^{2} {a_0} e f^{2} g {i_0} j k^{2} -a^{2} {b_0}^{2} f^{3} j^{2} k^{2}-2 a^{2} {b_0} d f^{2} {i_0} j^{2} k^{2}-2 a^{2} {b_0} e f^{2} g {i_0} j k^{2}+2 a {a_0}^{2} b e f^{2} i j k^{2} +2 a {a_0}^{2} e f^{2} h {i_0} j k r+4 a {a_0} b {b_0} e f^{2} i j k^{2}+4 a {a_0} b d e f i {i_0} j k^{2}+4 a {a_0} b e^{2} f g i {i_0} k^{2} +8 a {a_0} {b_0} e f^{2} h {i_0} j k r+4 a b {b_0} d e f i {i_0} j k^{2}+2 a {b_0}^{2} e f^{2} h {i_0} j k r-3 {a_0}^{2} b^{2} e^{2} f i^{2} k^{2} +3 {a_0}^{2} b^{2} e f^{2} h^{2} k^{2}-{a_0}^{2} d f^{2} i^{2} k^{2} r^{2}+3 {a_0}^{2} e^{3} i^{2} {i_0}^{2} r^{2}-6 {a_0} b^{2} d e^{2} i^{2} {i_0} k^{2} -4 {a_0} {b_0} d f^{2} i^{2} k^{2} r^{2}+6 {a_0} {b_0} e^{3} i^{2} {i_0}^{2} r^{2}-2 {a_0} d^{2} f i^{2} {i_0} k^{2} r^{2}-{b_0}^{2} d f^{2} i^{2} k^{2} r^{2} +{b_0}^{2} e^{3} i^{2} {i_0}^{2} r^{2}-4 {b_0} d^{2} f i^{2} {i_0} k^{2} r^{2}-d^{3} i^{2} {i_0}^{2} k^{2} r^{2}+2 a^{2} {a_0} f^{3} g j k^{2} +2 a^{2} {b_0} f^{3} g j k^{2}+2 a^{2} d f^{2} g {i_0} j k^{2}+a^{2} e f^{2} g^{2} {i_0} k^{2}-2 a {a_0}^{2} f^{3} h j k r -2 a {a_0} b d f^{2} i j k^{2}-4 a {a_0} b e f^{2} g i k^{2}-8 a {a_0} {b_0} f^{3} h j k r-2 a {a_0} d f^{2} h {i_0} j k r -4 a {a_0} e f^{2} g h {i_0} k r-2 a b {b_0} d f^{2} i j k^{2}-2 a b d^{2} f i {i_0} j k^{2}-4 a b d e f g i {i_0} k^{2} -2 a {b_0}^{2} f^{3} h j k r-4 a {b_0} d f^{2} h {i_0} j k r-2 a {b_0} e f^{2} g h {i_0} k r-6 {a_0}^{2} e^{2} f i^{2} {i_0} r^{2} +3 {a_0}^{2} e f^{2} h^{2} {i_0} r^{2}+4 {a_0} b^{2} d e f i^{2} k^{2}-2 {a_0} b^{2} d f^{2} h^{2} k^{2}-12 {a_0} {b_0} e^{2} f i^{2} {i_0} r^{2} +6 {a_0} {b_0} e f^{2} h^{2} {i_0} r^{2}-6 {a_0} d e^{2} i^{2} {i_0}^{2} r^{2}+3 b^{2} d^{2} e i^{2} {i_0} k^{2}-2 {b_0}^{2} e^{2} f i^{2} {i_0} r^{2} +{b_0}^{2} e f^{2} h^{2} {i_0} r^{2}-6 {b_0} d e^{2} i^{2} {i_0}^{2} r^{2}-a^{2} f^{3} g^{2} k^{2}+4 a {a_0} f^{3} g h k r+2 a b d f^{2} g i k^{2} +2 a {b_0} f^{3} g h k r+2 a d f^{2} g h {i_0} k r+3 {a_0}^{2} e f^{2} i^{2} r^{2}-3 {a_0}^{2} f^{3} h^{2} r^{2}+6 {a_0} {b_0} e f^{2} i^{2} r^{2} -6 {a_0} {b_0} f^{3} h^{2} r^{2}+8 {a_0} d e f i^{2} {i_0} r^{2}-2 {a_0} d f^{2} h^{2} {i_0} r^{2}-b^{2} d^{2} f i^{2} k^{2}+{b_0}^{2} e f^{2} i^{2} r^{2} -{b_0}^{2} f^{3} h^{2} r^{2}+8 {b_0} d e f i^{2} {i_0} r^{2}-2 {b_0} d f^{2} h^{2} {i_0} r^{2}+3 d^{2} e i^{2} {i_0}^{2} r^{2}-2 {a_0} d f^{2} i^{2} r^{2} -2 {b_0} d f^{2} i^{2} r^{2}-2 d^{2} f i^{2} {i_0} r^{2}$
    
    \item[$\mathbf{p_7} = $] $-4 a {a_0}^{2} {b_0} e f^{2} h {i_0} j k^{2} r-4 a {a_0} {b_0}^{2} e f^{2} h {i_0} j k^{2} r-2 {a_0}^{3} e^{3} i^{2} {i_0}^{2} k r^{2}-12 {a_0}^{2} {b_0} e^{3} i^{2} {i_0}^{2} k r^{2} -6 {a_0} {b_0}^{2} e^{3} i^{2} {i_0}^{2} k r^{2}+4 a {a_0}^{2} {b_0} f^{3} h j k^{2} r+2 a {a_0}^{2} e f^{2} g h {i_0} k^{2} r+4 a {a_0} {b_0}^{2} f^{3} h j k^{2} r +4 a {a_0} {b_0} d f^{2} h {i_0} j k^{2} r+4 a {a_0} {b_0} e f^{2} g h {i_0} k^{2} r+2 a {b_0}^{2} d f^{2} h {i_0} j k^{2} r+4 {a_0}^{3} e^{2} f i^{2} {i_0} k r^{2} -2 {a_0}^{3} e f^{2} h^{2} {i_0} k r^{2}+24 {a_0}^{2} {b_0} e^{2} f i^{2} {i_0} k r^{2}-12 {a_0}^{2} {b_0} e f^{2} h^{2} {i_0} k r^{2}+6 {a_0}^{2} d e^{2} i^{2} {i_0}^{2} k r^{2} +12 {a_0} {b_0}^{2} e^{2} f i^{2} {i_0} k r^{2}-6 {a_0} {b_0}^{2} e f^{2} h^{2} {i_0} k r^{2}+24 {a_0} {b_0} d e^{2} i^{2} {i_0}^{2} k r^{2}+6 {b_0}^{2} d e^{2} i^{2} {i_0}^{2} k r^{2} -2 a {a_0}^{2} f^{3} g h k^{2} r-4 a {a_0} {b_0} f^{3} g h k^{2} r-2 a {a_0} d f^{2} g h {i_0} k^{2} r-2 a {b_0} d f^{2} g h {i_0} k^{2} r -2 {a_0}^{3} e f^{2} i^{2} k r^{2}+2 {a_0}^{3} f^{3} h^{2} k r^{2}-12 {a_0}^{2} {b_0} e f^{2} i^{2} k r^{2}+12 {a_0}^{2} {b_0} f^{3} h^{2} k r^{2} -8 {a_0}^{2} d e f i^{2} {i_0} k r^{2}+2 {a_0}^{2} d f^{2} h^{2} {i_0} k r^{2}-6 {a_0} {b_0}^{2} e f^{2} i^{2} k r^{2}+6 {a_0} {b_0}^{2} f^{3} h^{2} k r^{2} -32 {a_0} {b_0} d e f i^{2} {i_0} k r^{2}+8 {a_0} {b_0} d f^{2} h^{2} {i_0} k r^{2}-6 {a_0} d^{2} e i^{2} {i_0}^{2} k r^{2}-8 {b_0}^{2} d e f i^{2} {i_0} k r^{2} +2 {b_0}^{2} d f^{2} h^{2} {i_0} k r^{2}-12 {b_0} d^{2} e i^{2} {i_0}^{2} k r^{2}+2 {a_0}^{2} d f^{2} i^{2} k r^{2}+8 {a_0} {b_0} d f^{2} i^{2} k r^{2} +4 {a_0} d^{2} f i^{2} {i_0} k r^{2}+2 {b_0}^{2} d f^{2} i^{2} k r^{2}+8 {b_0} d^{2} f i^{2} {i_0} k r^{2}+2 d^{3} i^{2} {i_0}^{2} k r^{2}$

    \item[$\mathbf{p_8} = $] $2 {a_0}^{3} {b_0} e^{3} i^{2} {i_0}^{2} k^{2} r^{2}+3 {a_0}^{2} {b_0}^{2} e^{3} i^{2} {i_0}^{2} k^{2} r^{2}-4 {a_0}^{3} {b_0} e^{2} f i^{2} {i_0} k^{2} r^{2}+2 {a_0}^{3} {b_0} e f^{2} h^{2} {i_0} k^{2} r^{2} -6 {a_0}^{2} {b_0}^{2} e^{2} f i^{2} {i_0} k^{2} r^{2}+3 {a_0}^{2} {b_0}^{2} e f^{2} h^{2} {i_0} k^{2} r^{2}-6 {a_0}^{2} {b_0} d e^{2} i^{2} {i_0}^{2} k^{2} r^{2} -6 {a_0} {b_0}^{2} d e^{2} i^{2} {i_0}^{2} k^{2} r^{2}+a^{2} {a_0} {b_0}^{2} e f^{2} {i_0} j^{2} k^{2}-2 a {a_0}^{2} b {b_0} e^{2} f i {i_0} j k^{2}+{a_0}^{3} b^{2} e^{3} i^{2} {i_0} k^{2} +2 {a_0}^{3} {b_0} e f^{2} i^{2} k^{2} r^{2}-2 {a_0}^{3} {b_0} f^{3} h^{2} k^{2} r^{2}+3 {a_0}^{2} {b_0}^{2} e f^{2} i^{2} k^{2} r^{2}-3 {a_0}^{2} {b_0}^{2} f^{3} h^{2} k^{2} r^{2} +8 {a_0}^{2} {b_0} d e f i^{2} {i_0} k^{2} r^{2}-2 {a_0}^{2} {b_0} d f^{2} h^{2} {i_0} k^{2} r^{2}+8 {a_0} {b_0}^{2} d e f i^{2} {i_0} k^{2} r^{2} -2 {a_0} {b_0}^{2} d f^{2} h^{2} {i_0} k^{2} r^{2}+6 {a_0} {b_0} d^{2} e i^{2} {i_0}^{2} k^{2} r^{2}+3 {b_0}^{2} d^{2} e i^{2} {i_0}^{2} k^{2} r^{2}-a^{2} {a_0} {b_0}^{2} f^{3} j^{2} k^{2} -2 a^{2} {a_0} {b_0} e f^{2} g {i_0} j k^{2}-a^{2} {b_0}^{2} d f^{2} {i_0} j^{2} k^{2}+2 a {a_0}^{2} b {b_0} e f^{2} i j k^{2}+2 a {a_0}^{2} b e^{2} f g i {i_0} k^{2} +4 a {a_0}^{2} {b_0} e f^{2} h {i_0} j k r+4 a {a_0} b {b_0} d e f i {i_0} j k^{2}+4 a {a_0} {b_0}^{2} e f^{2} h {i_0} j k r-{a_0}^{3} b^{2} e^{2} f i^{2} k^{2} +{a_0}^{3} b^{2} e f^{2} h^{2} k^{2}+{a_0}^{3} e^{3} i^{2} {i_0}^{2} r^{2}-3 {a_0}^{2} b^{2} d e^{2} i^{2} {i_0} k^{2}-2 {a_0}^{2} {b_0} d f^{2} i^{2} k^{2} r^{2} +6 {a_0}^{2} {b_0} e^{3} i^{2} {i_0}^{2} r^{2}-2 {a_0} {b_0}^{2} d f^{2} i^{2} k^{2} r^{2}+3 {a_0} {b_0}^{2} e^{3} i^{2} {i_0}^{2} r^{2}-4 {a_0} {b_0} d^{2} f i^{2} {i_0} k^{2} r^{2} -2 {b_0}^{2} d^{2} f i^{2} {i_0} k^{2} r^{2}-2 {b_0} d^{3} i^{2} {i_0}^{2} k^{2} r^{2}+2 a^{2} {a_0} {b_0} f^{3} g j k^{2}+a^{2} {a_0} e f^{2} g^{2} {i_0} k^{2} +2 a^{2} {b_0} d f^{2} g {i_0} j k^{2}-2 a {a_0}^{2} b e f^{2} g i k^{2}-4 a {a_0}^{2} {b_0} f^{3} h j k r-2 a {a_0}^{2} e f^{2} g h {i_0} k r -2 a {a_0} b {b_0} d f^{2} i j k^{2}-4 a {a_0} b d e f g i {i_0} k^{2}-4 a {a_0} {b_0}^{2} f^{3} h j k r-4 a {a_0} {b_0} d f^{2} h {i_0} j k r -4 a {a_0} {b_0} e f^{2} g h {i_0} k r-2 a b {b_0} d^{2} f i {i_0} j k^{2}-2 a {b_0}^{2} d f^{2} h {i_0} j k r-2 {a_0}^{3} e^{2} f i^{2} {i_0} r^{2} +{a_0}^{3} e f^{2} h^{2} {i_0} r^{2}+2 {a_0}^{2} b^{2} d e f i^{2} k^{2}-{a_0}^{2} b^{2} d f^{2} h^{2} k^{2}-12 {a_0}^{2} {b_0} e^{2} f i^{2} {i_0} r^{2} +6 {a_0}^{2} {b_0} e f^{2} h^{2} {i_0} r^{2}-3 {a_0}^{2} d e^{2} i^{2} {i_0}^{2} r^{2}+3 {a_0} b^{2} d^{2} e i^{2} {i_0} k^{2}-6 {a_0} {b_0}^{2} e^{2} f i^{2} {i_0} r^{2} +3 {a_0} {b_0}^{2} e f^{2} h^{2} {i_0} r^{2}-12 {a_0} {b_0} d e^{2} i^{2} {i_0}^{2} r^{2}-3 {b_0}^{2} d e^{2} i^{2} {i_0}^{2} r^{2}-a^{2} {a_0} f^{3} g^{2} k^{2} -a^{2} d f^{2} g^{2} {i_0} k^{2}+2 a {a_0}^{2} f^{3} g h k r+2 a {a_0} b d f^{2} g i k^{2}+4 a {a_0} {b_0} f^{3} g h k r+2 a {a_0} d f^{2} g h {i_0} k r +2 a b d^{2} f g i {i_0} k^{2}+2 a {b_0} d f^{2} g h {i_0} k r+{a_0}^{3} e f^{2} i^{2} r^{2}-{a_0}^{3} f^{3} h^{2} r^{2}+6 {a_0}^{2} {b_0} e f^{2} i^{2} r^{2} -6 {a_0}^{2} {b_0} f^{3} h^{2} r^{2}+4 {a_0}^{2} d e f i^{2} {i_0} r^{2}-{a_0}^{2} d f^{2} h^{2} {i_0} r^{2}-{a_0} b^{2} d^{2} f i^{2} k^{2}+3 {a_0} {b_0}^{2} e f^{2} i^{2} r^{2} -3 {a_0} {b_0}^{2} f^{3} h^{2} r^{2}+16 {a_0} {b_0} d e f i^{2} {i_0} r^{2}-4 {a_0} {b_0} d f^{2} h^{2} {i_0} r^{2}+3 {a_0} d^{2} e i^{2} {i_0}^{2} r^{2} -b^{2} d^{3} i^{2} {i_0} k^{2}+4 {b_0}^{2} d e f i^{2} {i_0} r^{2}-{b_0}^{2} d f^{2} h^{2} {i_0} r^{2}+6 {b_0} d^{2} e i^{2} {i_0}^{2} r^{2}-{a_0}^{2} d f^{2} i^{2} r^{2} -4 {a_0} {b_0} d f^{2} i^{2} r^{2}-2 {a_0} d^{2} f i^{2} {i_0} r^{2}-{b_0}^{2} d f^{2} i^{2} r^{2}-4 {b_0} d^{2} f i^{2} {i_0} r^{2}-d^{3} i^{2} {i_0}^{2} r^{2}$

    \item[$\mathbf{p_9} = $] $-2 a {a_0}^{2} {b_0}^{2} e f^{2} h {i_0} j k^{2} r-4 {a_0}^{3} {b_0} e^{3} i^{2} {i_0}^{2} k r^{2}-6 {a_0}^{2} {b_0}^{2} e^{3} i^{2} {i_0}^{2} k r^{2} +2 a {a_0}^{2} {b_0}^{2} f^{3} h j k^{2} r+2 a {a_0}^{2} {b_0} e f^{2} g h {i_0} k^{2} r+2 a {a_0} {b_0}^{2} d f^{2} h {i_0} j k^{2} r +8 {a_0}^{3} {b_0} e^{2} f i^{2} {i_0} k r^{2}-4 {a_0}^{3} {b_0} e f^{2} h^{2} {i_0} k r^{2}+12 {a_0}^{2} {b_0}^{2} e^{2} f i^{2} {i_0} k r^{2} -6 {a_0}^{2} {b_0}^{2} e f^{2} h^{2} {i_0} k r^{2}+12 {a_0}^{2} {b_0} d e^{2} i^{2} {i_0}^{2} k r^{2}+12 {a_0} {b_0}^{2} d e^{2} i^{2} {i_0}^{2} k r^{2}  -2 a {a_0}^{2} {b_0} f^{3} g h k^{2} r-2 a {a_0} {b_0} d f^{2} g h {i_0} k^{2} r-4 {a_0}^{3} {b_0} e f^{2} i^{2} k r^{2}+4 {a_0}^{3} {b_0} f^{3} h^{2} k r^{2} -6 {a_0}^{2} {b_0}^{2} e f^{2} i^{2} k r^{2}+6 {a_0}^{2} {b_0}^{2} f^{3} h^{2} k r^{2}-16 {a_0}^{2} {b_0} d e f i^{2} {i_0} k r^{2}+4 {a_0}^{2} {b_0} d f^{2} h^{2} {i_0} k r^{2} -16 {a_0} {b_0}^{2} d e f i^{2} {i_0} k r^{2}+4 {a_0} {b_0}^{2} d f^{2} h^{2} {i_0} k r^{2}-12 {a_0} {b_0} d^{2} e i^{2} {i_0}^{2} k r^{2}-6 {b_0}^{2} d^{2} e i^{2} {i_0}^{2} k r^{2} +4 {a_0}^{2} {b_0} d f^{2} i^{2} k r^{2}+4 {a_0} {b_0}^{2} d f^{2} i^{2} k r^{2}+8 {a_0} {b_0} d^{2} f i^{2} {i_0} k r^{2}+4 {b_0}^{2} d^{2} f i^{2} {i_0} k r^{2} +4 {b_0} d^{3} i^{2} {i_0}^{2} k r^{2}$

    \item[$\mathbf{p_{10}} = $] ${a_0}^{3} {b_0}^{2} e^{3} i^{2} {i_0}^{2} k^{2} r^{2}-2 {a_0}^{3} {b_0}^{2} e^{2} f i^{2} {i_0} k^{2} r^{2}+{a_0}^{3} {b_0}^{2} e f^{2} h^{2} {i_0} k^{2} r^{2}  -3 {a_0}^{2} {b_0}^{2} d e^{2} i^{2} {i_0}^{2} k^{2} r^{2}+{a_0}^{3} {b_0}^{2} e f^{2} i^{2} k^{2} r^{2}-{a_0}^{3} {b_0}^{2} f^{3} h^{2} k^{2} r^{2}+4 {a_0}^{2} {b_0}^{2} d e f i^{2} {i_0} k^{2} r^{2} -{a_0}^{2} {b_0}^{2} d f^{2} h^{2} {i_0} k^{2} r^{2}+3 {a_0} {b_0}^{2} d^{2} e i^{2} {i_0}^{2} k^{2} r^{2}+2 a {a_0}^{2} {b_0}^{2} e f^{2} h {i_0} j k r+2 {a_0}^{3} {b_0} e^{3} i^{2} {i_0}^{2} r^{2} -{a_0}^{2} {b_0}^{2} d f^{2} i^{2} k^{2} r^{2}+3 {a_0}^{2} {b_0}^{2} e^{3} i^{2} {i_0}^{2} r^{2}-2 {a_0} {b_0}^{2} d^{2} f i^{2} {i_0} k^{2} r^{2}-{b_0}^{2} d^{3} i^{2} {i_0}^{2} k^{2} r^{2} -2 a {a_0}^{2} {b_0}^{2} f^{3} h j k r-2 a {a_0}^{2} {b_0} e f^{2} g h {i_0} k r-2 a {a_0} {b_0}^{2} d f^{2} h {i_0} j k r-4 {a_0}^{3} {b_0} e^{2} f i^{2} {i_0} r^{2} +2 {a_0}^{3} {b_0} e f^{2} h^{2} {i_0} r^{2}-6 {a_0}^{2} {b_0}^{2} e^{2} f i^{2} {i_0} r^{2}+3 {a_0}^{2} {b_0}^{2} e f^{2} h^{2} {i_0} r^{2}-6 {a_0}^{2} {b_0} d e^{2} i^{2} {i_0}^{2} r^{2} -6 {a_0} {b_0}^{2} d e^{2} i^{2} {i_0}^{2} r^{2}+2 a {a_0}^{2} {b_0} f^{3} g h k r+2 a {a_0} {b_0} d f^{2} g h {i_0} k r+2 {a_0}^{3} {b_0} e f^{2} i^{2} r^{2} -2 {a_0}^{3} {b_0} f^{3} h^{2} r^{2}+3 {a_0}^{2} {b_0}^{2} e f^{2} i^{2} r^{2}-3 {a_0}^{2} {b_0}^{2} f^{3} h^{2} r^{2}+8 {a_0}^{2} {b_0} d e f i^{2} {i_0} r^{2} -2 {a_0}^{2} {b_0} d f^{2} h^{2} {i_0} r^{2}+8 {a_0} {b_0}^{2} d e f i^{2} {i_0} r^{2}-2 {a_0} {b_0}^{2} d f^{2} h^{2} {i_0} r^{2}+6 {a_0} {b_0} d^{2} e i^{2} {i_0}^{2} r^{2} +3 {b_0}^{2} d^{2} e i^{2} {i_0}^{2} r^{2}-2 {a_0}^{2} {b_0} d f^{2} i^{2} r^{2}-2 {a_0} {b_0}^{2} d f^{2} i^{2} r^{2}-4 {a_0} {b_0} d^{2} f i^{2} {i_0} r^{2} -2 {b_0}^{2} d^{2} f i^{2} {i_0} r^{2}-2 {b_0} d^{3} i^{2} {i_0}^{2} r^{2}$

    \item[$\mathbf{p_{11}} = $] $-2 {a_0}^{3} {b_0}^{2} e^{3} i^{2} {i_0}^{2} k r^{2}+4 {a_0}^{3} {b_0}^{2} e^{2} f i^{2} {i_0} k r^{2}-2 {a_0}^{3} {b_0}^{2} e f^{2} h^{2} {i_0} k r^{2} +6 {a_0}^{2} {b_0}^{2} d e^{2} i^{2} {i_0}^{2} k r^{2}-2 {a_0}^{3} {b_0}^{2} e f^{2} i^{2} k r^{2}+2 {a_0}^{3} {b_0}^{2} f^{3} h^{2} k r^{2} -8 {a_0}^{2} {b_0}^{2} d e f i^{2} {i_0} k r^{2} +2 {a_0}^{2} {b_0}^{2} d f^{2} h^{2} {i_0} k r^{2}-6 {a_0} {b_0}^{2} d^{2} e i^{2} {i_0}^{2} k r^{2} +2 {a_0}^{2} {b_0}^{2} d f^{2} i^{2} k r^{2}+4 {a_0} {b_0}^{2} d^{2} f i^{2} {i_0} k r^{2} +2 {b_0}^{2} d^{3} i^{2} {i_0}^{2} k r^{2}$

    \item[$\mathbf{p_{12}} = $] ${a_0}^{3} {b_0}^{2} e^{3} i^{2} {i_0}^{2} r^{2}-2 {a_0}^{3} {b_0}^{2} e^{2} f i^{2} {i_0} r^{2} +{a_0}^{3} {b_0}^{2} e f^{2} h^{2} {i_0} r^{2}-3 {a_0}^{2} {b_0}^{2} d e^{2} i^{2} {i_0}^{2} r^{2} +{a_0}^{3} {b_0}^{2} e f^{2} i^{2} r^{2}-{a_0}^{3} {b_0}^{2} f^{3} h^{2} r^{2}+4 {a_0}^{2} {b_0}^{2} d e f i^{2} {i_0} r^{2}-{a_0}^{2} {b_0}^{2} d f^{2} h^{2} {i_0} r^{2} +3 {a_0} {b_0}^{2} d^{2} e i^{2} {i_0}^{2} r^{2}-{a_0}^{2} {b_0}^{2} d f^{2} i^{2} r^{2}-2 {a_0} {b_0}^{2} d^{2} f i^{2} {i_0} r^{2}-{b_0}^{2} d^{3} i^{2} {i_0}^{2} r^{2}$
\end{enumerate}}

\end{document}